\documentclass{amsart}
\usepackage{amssymb,amsmath,amsthm}
\newtheorem{theorem}{Theorem}[section]
\newtheorem{lemma}[theorem]{Lemma}
\newtheorem{remark}{Remark}[section]
\usepackage{esint}
\usepackage{xcolor}

\usepackage{graphicx}

\usepackage{color}
\usepackage{hyperref}
\hypersetup{
	colorlinks = true,%
	citecolor = blue,
	filecolor=red,%
	linkcolor = [rgb]{0.65,0.0,0.0},%
	anchorcolor = red,
	pagecolor = red,
	urlcolor= [rgb]{0.65,0.0,0.0}
	linktocpage=true,
	pdfpagelabels=true,
	bookmarksnumbered=true,
}

\numberwithin{equation}{section}
\allowdisplaybreaks[1]

\begin{document}

\title{Elliptic problem involving finite many critical exponents in $\mathbb{R}^{N}$}




\author{Yu Su}
\address{School of Mathematics and Statistics, Central South University, Changsha, 410083 Hunan, P.R.China.}
\curraddr{}
\email{yizai52@qq.com}
\thanks{}

\author{Haibo Chen}
\address{School of Mathematics and Statistics, Central South University, Changsha, 410083 Hunan, P.R.China.}
\curraddr{}
\email{math\_chb@163.com}
\thanks{This research was supported by National Natural Science Foundation of China 11671403.}



\date{}

\dedicatory{}

\begin{abstract}
In this paper,
we consider  the following problem
$$
-\Delta u
-\zeta
\frac{u}{|x|^{2}}
=
\sum_{i=1}^{k}
\left(
\int_{\mathbb{R}^{N}}
\frac{|u|^{2^{*}_{\alpha_{i}}}}{|x-y|^{\alpha_{i}}}
\mathrm{d}y
\right)
|u|^{2^{*}_{\alpha_{i}}-2}u
+
|u|^{2^{*}-2}u
,
\mathrm{~in~}
\mathbb{R}^{N},
$$
where
$N\geqslant3$,
$\zeta\in(0,\frac{(N-2)^{2}}{4})$,
$2^{*}=\frac{2N}{N-2}$
is the critical Sobolev exponent,
and
$2^{*}_{\alpha_{i}}=\frac{2N-\alpha_{i}}{N-2}$
($i=1,\ldots,k$)
are the  Hardy--Littlewood--Sobolev critical upper exponents.
The parameters
$\alpha_{i}$
($i=1,\ldots,k$)
satisfy some suitable assumptions.
By using Coulomb--Sobolev space,
endpoint refined Sobolev inequality
and variational methods,
we
establish the existence of nontrivial solutions.
Our result extends the ones in Yang and Wu [Adv. Nonlinear Stud. (2017) \cite{Yang2017}].
\par
\end{abstract}

\maketitle

\section{Introduction}
In this paper,
we  consider the following problem:
$$
-\Delta u
-\zeta
\frac{u}{|x|^{2}}
=
\sum_{i=1}^{k}
\left(
\int_{\mathbb{R}^{N}}
\frac{|u|^{2^{*}_{\alpha_{i}}}}{|x-y|^{\alpha_{i}}}
\mathrm{d}y
\right)
|u|^{2^{*}_{\alpha_{i}}-2}u
+
|u|^{2^{*}-2}u
,
\mathrm{~in~}
\mathbb{R}^{N},
\eqno(\mathcal{P})
$$
where
$N\geqslant3$,
$\zeta\in(0,\frac{(N-2)^{2}}{4})$,
$2^{*}=
\frac{2N}{N-2}$
is the critical Sobolev  exponent,
$2^{*}_{\alpha_{i}}=
\frac{2N-\alpha_{i}}{N-2}$
($i=1,\ldots,k$)
are the Hardy--Littlewood--Sobolev critical upper  exponents,
and
the parameters $\alpha_{i}$
($i=1,\ldots,k$) satisfy the following assumptions:

\begin{itemize}
\item[{\bf ($H_{1}$)}]
$0<\alpha_{1}<\alpha_{2}<\cdots<\alpha_{k}<N$
($k\in \mathbb{N}$,
$2\leqslant k<\infty$);

\item[{\bf ($H_{2}$)}]
$
\frac{
\pi N (N-2)
\left(
\frac{\Gamma(\frac{N}{2})}{\Gamma(N)}
\right)^{\frac{2}{N}}
}
{
\left(
\pi^{\frac{\alpha_{i}}{2}}
\frac{\Gamma(\frac{N}{2}-\frac{\alpha_{i}}{2})}{\Gamma(N-\frac{\alpha_{i}}{2})}
\left(
\frac{\Gamma(\frac{N}{2})}{\Gamma(N)}
\right)^{\frac{\alpha_{i}-N}{N}}
\right)
^{\frac{1}{2^{*}_{\alpha_{i}}}}}
\geqslant1$,
for all $i=1,\ldots,k$.
\end{itemize}

According to
\eqref{14}
and
\eqref{15},
we could see that the assumption $(H_{2})$
is equivalent to $\tilde{S}_{\alpha_{i}}\geqslant1$ for all $i=1,\ldots,k$.
For any $N\geqslant3$,
the value of
$\tilde{S}_{\alpha}$
is dependent on the parameters
$N$
and
$\alpha$ (see Fig 1.).
\begin{center}
\includegraphics[height=5cm]{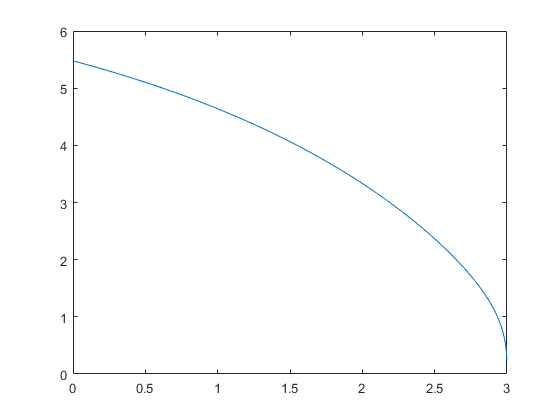}\\
\centerline{Fig 1. The $x$--axis is $\alpha$, and the $y$--axis is $\tilde{S}_{\alpha}$, where $N=3$ and $\alpha\in(0,N)$. }
\end{center}
From Fig 1,
there exists a constant $0<\bar{\alpha}<3$ such that $\tilde{S}_{\bar{\alpha}}=1$.
For any $\alpha\in (0,\bar{\alpha}]$,
we have $1\leqslant\tilde{S}_{\alpha}<\infty$.
For any
$\alpha\in (\bar{\alpha},3)$,
we get $0<\tilde{S}_{\alpha}<1$.
In this paper,
we just study the case of $1\leqslant\tilde{S}_{\alpha}<\infty$.

Problem
$(\mathcal{P})$
is related to the nonlinear Choquard equation as follows:
\begin{equation}\label{1}
-\Delta u
+
V(x)u
=
\left(
|x|^{\alpha}*|u|^{q}
\right)
|u|^{q-2}u,
\mathrm{~in~}
\mathbb{R}^{N},
\end{equation}
where
$\frac{2N-\alpha}{N}\leqslant q\leqslant\frac{2N-\alpha}{N-2}$
and
$\alpha\in(0,N)$.
For
$q=2$
and
$\alpha=1$,
the problem
(\ref{1})
goes back to the description of the quantum theory of a polaron at rest by Pekar in 1954
\cite{Pekar1954}
and the modeling of an electron trapped in its own hole in
1976
in the work of
Choquard,
as a certain approximation to Hartree--Fock theory of one--component plasma
\cite{Penrose1996}.
The existence and qualitative properties of solutions of Choquard type equations
(\ref{1})
have been widely studied in the last decades
(see \cite{Moroz2016}).

For Laplacian with nonlocal Hartree--type nonlinearities,
the problem has attracted a lot of interest.
Gao and Yang \cite{Gao2016} investigated the following critical Choquard equation:
\begin{equation}\label{2}
\begin{aligned}
-\Delta u
=
\left(
\int_{\mathbb{R}^{N}}
\frac{|u|^{2^{*}_{\alpha}}}{|x-y|^{\alpha}}
\mathrm{d}y
\right)
|u|^{2^{*}_{\alpha}-2}u
+\lambda u,
\mathrm{~in~}
\Omega,
\end{aligned}
\end{equation}
where
$\Omega$
is a bounded domain of
$\mathbb{R}^{N}$,
with lipschitz boundary,
$N\geqslant3$,
$\alpha\in(0,N)$
and
$\lambda>0$.
By using variational methods,
they established the existence, multiplicity and nonexistence of nontrivial solutions to equation (\ref{2}).
Alves, Gao,  Squassina and Yang
\cite{O.Alves2017}
studied the following
singularly perturbed critical Choquard equation:
\begin{equation*}
\begin{aligned}
-\varepsilon^{2}\Delta u
+
V(x)u
=
\varepsilon^{\alpha-3}
\left(
\int_{\mathbb{R}^{N}}
\frac{Q(y)G(u(y))}{|x-y|^{\alpha}}
\mathrm{d}y
\right)
Q(x)g(u),
\mathrm{~in~}
\mathbb{R}^{3},
\end{aligned}
\end{equation*}
where
$0<\alpha<3$,
$\varepsilon$
is a positive parameter,
$V,Q$
are two continuous real functions on
$\mathbb{R}^{3}$
and
$G$
is the primitive of
$g$
which is of critical growth due to the Hardy--Littlewood--Sobolev inequality.
Under suitable assumptions on
$g$,
they first establish the existence of ground states for the critical Choquard equation with constant coefficient.
They also established existence and multiplicity of semi--classical solutions and characterize the concentration behavior by variational methods.
For details and recent works,
we refer to \cite{Gao2017JMAA,Mercuri2016} and the references therein.

For fractional Laplacian with nonlocal Hartree--type nonlinearities,
the problem has attracted a lot of interest.
D'Avenia,  Siciliano and Squassina
\cite{d'Avenia2015}
considered the following fractional Choquard equation:
\begin{equation}\label{5}
\begin{aligned}
(-\Delta)^{s} u
+
\omega u
=
\left(
\mathcal{K}_{\alpha}\ast|u|^{q}
\right)
|u|^{q-2}u,
\mathrm{~in~}
\mathbb{R}^{N},
\end{aligned}
\end{equation}
where
$N\geqslant3$,
$s\in(0,1)$,
$\omega\geqslant0$,
$\alpha\in(0,N)$
and
$q\in(\frac{2N-\alpha}{N},\frac{2N-\alpha}{N-2s})$.
In particularly,
if
$\omega=0$,
$\alpha=4s$
and
$q=2$,
then peoblem
(\ref{5})
becomes a fractional Choquard euqation with upper critical exponent in the sense of Hardy--Littlewood--Sobolev inequatlity as follows:
\begin{equation}\label{6}
\begin{aligned}
(-\Delta)^{s} u
=
\left(
\int_{\mathbb{R}^{N}}
\frac{|u|^{2}}{|x-y|^{4s}}
\mathrm{d}y
\right)
u,
\mathrm{~in~}
\mathbb{R}^{N}.
\end{aligned}
\end{equation}
D'Avenia,  Siciliano and Squassina in
\cite{d'Avenia2015}
obtained regularity, existence, nonexistence of nontrivial solutions to problem (\ref{5}) and problem (\ref{6}).
Mukherjee and Sreenadh \cite{Mukherjee2017Fractional} investigated the following fractional Choquard equation:
\begin{equation}\label{4}
\begin{aligned}
(-\Delta)^{s} u
=
\left(
\int_{\mathbb{R}^{N}}
\frac{|u|^{2^{*}_{h,\alpha}}}{|x-y|^{\alpha}}
\mathrm{d}y
\right)
|u|^{2^{*}_{h,\alpha}-2}u
+\lambda u,
\mathrm{~in~}
\Omega,
\end{aligned}
\end{equation}
where
$\Omega$
is a bounded domain of
$\mathbb{R}^{N}$
with $C^{1,1}$ boundary,
$N\geqslant3$,
$s\in(0,1)$,
$\alpha\in(0,N)$,
$\lambda>0$
and
$2^{*}_{h,\alpha}=\frac{2N-\alpha}{N-2s}$.
Applying variational methods,
they established the existence, multiplicity and nonexistence of nontrivial solutions to problem (\ref{4}).

Recently,
Yang and Wu \cite{Yang2017}
studied the following nonlocal elliptic problems:
\begin{equation}\label{7}
\begin{aligned}
(-\Delta)^{s} u
-\frac{\zeta u}{|x|^{2s}}
=
\left(
\int_{\mathbb{R}^{N}}
\frac{|u|^{2^{*}_{h,\alpha}}}{|x-y|^{\alpha}}
\mathrm{d}y
\right)
|u|^{2^{*}_{h,\alpha}-2}u
+
\left(
\int_{\mathbb{R}^{N}}
\frac{|u|^{2^{*}_{h,\beta}}}{|x-y|^{\beta}}
\mathrm{d}y
\right)
|u|^{2^{*}_{h,\beta}-2}u
,
\mathrm{~in~}
\mathbb{R}^{N},
\end{aligned}
\end{equation}
and
\begin{equation}\label{8}
\begin{aligned}
(-\Delta)^{s} u
-\frac{\zeta u}{|x|^{2s}}
=
\left(
\int_{\mathbb{R}^{N}}
\frac{|u|^{2^{*}_{h,\alpha}}}{|x-y|^{\alpha}}
\mathrm{d}y
\right)
|u|^{2^{*}_{h,\alpha}-2}u
+
\frac{|u|^{2^{*}_{s,\theta}-2}u}{|x|^{\theta}},
\mathrm{~in~}
\mathbb{R}^{N},
\end{aligned}
\end{equation}
where
$N\geqslant3$,
$s\in(0,1)$,
$\zeta\in
\left[
0,4^{s}\frac{\Gamma(\frac{N+2s}{4})}{\Gamma(\frac{N-2s}{4})}
\right)$,
$\alpha,\beta\in(N-2s,N)$,
$\theta\in(0,2s)$,
$2^{*}_{h,\alpha}=\frac{2N-\alpha}{N-2s}$
and
$2^{*}_{s,\theta}=\frac{2(N-\theta)}{N-2s}$.
By using a refinement of the Sobolev inequality which is related to
the Morrey space,
they showed the existence of nontrivial solutions for problem
(\ref{7})
and
problem
(\ref{8}).

In \cite{ZhangBL2017},
Wang, Zhang and Zhang
 extended the study of problem (\ref{8}) to the fractional Laplacian system as follows:
\begin{equation}\label{9}
\begin{aligned}
\begin{cases}
(-\Delta)^{s} u
-\!
\frac{\zeta u}{|x|^{2s}}
\!=\!
\left(
\int\limits_{\mathbb{R}^{N}}
\frac{|u|^{2^{*}_{h,\alpha}}}{|x-y|^{\alpha}}
\mathrm{d}y
\right)\!
|u|^{2^{*}_{h,\alpha}-2}u
\!+\!
\frac{|u|^{2^{*}_{s,\theta}-2}u}{|x|^{\theta}}
\!+\!
\frac{\beta\eta}{2^{*}_{s,\theta}}
\frac{
|u|^{\beta-2}u
|v|^{\gamma}}
{|x|^{\theta}}
,\\
(-\Delta)^{s} v
-\!
\frac{\zeta v}{|x|^{2s}}
\!=\!
\left(
\int\limits_{\mathbb{R}^{N}}
\frac{|v|^{2^{*}_{h,\alpha}}}{|x-y|^{\alpha}}
\mathrm{d}y
\right)\!
|v|^{2^{*}_{h,\alpha}-2}v
+\!
\frac{|v|^{2^{*}_{s,\theta}-2}v}{|x|^{\theta}}
+\!
\frac{\gamma\eta}{2^{*}_{s,\theta}}
\frac{|u|^{\beta}
|v|^{\gamma-2}v}
{|x|^{\theta}}
,
\end{cases}
\end{aligned}
\end{equation}
where
$N\geqslant3$,
$s\in(0,1)$,
$\zeta\in
\left[
0,4^{s}\frac{\Gamma(\frac{N+2s}{4})}{\Gamma(\frac{N-2s}{4})}
\right)$,
$\eta\in \mathbb{R}^{+}_{0}$,
$\alpha\in(N-2s,N)$,
$\theta\in(0,2s)$,
$\beta>1,\gamma>1$
and
$\beta+\gamma=2^{*}_{s,\theta}$.
By using variational methods,
they investigated the extremals of the corresponding best fractional Hardy--Sobolev constant and established the existence of solutions to problem (\ref{9}).

Moreover,
there are many other kinds of problems involving two critical nonlinearities,
such as the Laplacian
$-\Delta$
(see \cite{Li2012,Seok2018,Zhong2016}),
the $p$--Laplacian
$-\Delta_{p}$
(see \cite{Pucci2009}),
the biharmonic operator
$\Delta^{2}$ (see \cite{Bhakta2015}),
and the fractional operator
$(-\Delta)^{s}$
(see \cite{Ghoussoub2016}).
{\bf A natural and interesting question is:
For $s=1$,
can we extend the study of problem
(\ref{7})
to problem $(\mathcal{P})$?}
In this paper,
we give a positive answer to the question.
We need the following
inequalities.
\begin{lemma}\label{lemma1}
$\left.\right.$\cite[Hardy-Littlewood-Sobolev~inequality]{Lieb2001}
Let
$t,r>1$
and
$0<\alpha<N$
with
$\frac{1}{t}+\frac{1}{r}+\frac{\alpha}{N}=2$,
$f\in L^{t}(\mathbb{R}^{N})$
and
$h\in L^{r}(\mathbb{R}^{N})$.
There exists a sharp constant
$C(N,\alpha,t,r)>0$,
independent of
$f,g$
such that
$$\int_{\mathbb{R}^{N}}\int_{\mathbb{R}^{N}}
\frac{|f(x)||h(y)|}
{|x-y|^{\alpha}}
\mathrm{d}x\mathrm{d}y
\leqslant
C(N,\alpha,t,r)
\|f\|_{t}
\|h\|_{r}.$$
If
$t=r=\frac{2N}{2N-\alpha}$,
then
$$C(N,\alpha,t,r)
=
C(N,\alpha)
=
\pi^{\frac{\alpha}{2}}
\frac{\Gamma(\frac{N}{2}-\frac{\alpha}{2})}{\Gamma(N-\frac{\alpha}{2})}
\left\{
\frac{\Gamma(\frac{N}{2})}{\Gamma(N)}
\right\}^{\frac{\alpha-N}{N}}.$$
\end{lemma}
\begin{lemma}\label{lemma2}
$\left.\right.$\cite[Endpoint refined Sobolev inequality]{Mercuri2016}
Let
$\alpha\in(0,N)$.
Then there exists a  constant
$C_{1}>0$
such that the inequality
$$
\|u\|_{L^{2^{*}}(\mathbb{R}^{N})}
\leqslant
C_{1}
\|u\|_{D}^{\frac{(N-\alpha)(N-2)}{N(N+2-\alpha)}}
\left(
\int_{\mathbb{R}^{N}}
\int_{\mathbb{R}^{N}}
\frac{|u(x)|^{2^{*}_{\alpha}}|u(y)|^{2^{*}_{\alpha}}}{|x-y|^{\alpha}}
\mathrm{d}x
\mathrm{d}y
\right)
^{\frac{N-2}{N(N+2-\alpha)}},$$
holds for all
$u\in \mathcal{E}^{1,\alpha,2^{*}_{\alpha}}(\mathbb{R}^{N})$.
\end{lemma}
For the Coulomb--Sobolev space and endpoint refined Sobolev inequality,
there are two papers until now.
For Laplacian operator,
Mercuri, Moroz and Schaftingen \cite{Mercuri2016} introduced the Coulomb--Sobolev space and a family of associated optimal interpolation inequalities
(endpoint refined Sobolev inequality).
They established the existence of solutions of the nonlocal Schr\"{o}dinger--Poisson--Slater type equation in \cite{Mercuri2016}.
For fractional Laplacian operator,
Bellazzini, Ghimenti, Mercuri, Moroz and Schaftingen \cite{Bellazzini2016} studied the fractional Coulomb--Sobolev space and endpoint refined Sobolev inequality.

In this paper,
we apply Coulomb--Sobolev space and endpoint refined Sobolev to study problem $(\mathcal{P})$.
The main
result of this paper reads as follows.
\begin{theorem}\label{theorem1}
Let
$N\geqslant3$,
$(H_{1})$
and
$(H_{2})$
hold.
Then
problem
$(\mathcal{P})$
has a nonnegative solution
$\tilde{v}(x)$
for
$$\zeta
\in
\left(
\frac{(N-2)^{2}}{4}
\left(
1-
\frac{1}
{
(k+1)^{\frac{N-2}{N-1}}
\tilde{S}^{\frac{N}{N-1}}
}
\right)
,
\frac{(N-2)^{2}}{4}
\right).$$
Moreover,
set
$
\tilde{\tilde{v}}(x)
=
\frac{1}{|x|^{N-2}}
\tilde{v}
\left(
\frac{x}{|x|^{2}}
\right)$.
Then
$\tilde{\tilde{v}}(x)$
is a nonnegative solution of the problem
\begin{equation*}
\begin{aligned}
-
\Delta
\tilde{\tilde{v}}
-\zeta
\frac{\tilde{\tilde{v}}}{|x|^{2}}
=&
\sum\limits^{k}_{i=1}
\left(
\int_{\mathbb{R}^{N}}
\frac{|\tilde{\tilde{v}}|^{2^{*}_{\alpha_{i}}}}
{|x-y|^{\alpha_{i}}}
\mathrm{d}y
\right)
\left|
\tilde{\tilde{v}}
\right|^{2^{*}_{\alpha_{i}}-2}
\tilde{\tilde{v}}
+
\left|
\tilde{\tilde{v}}
\right|^{2^{*}-2}
\tilde{\tilde{v}},
~\mathrm{in}~\mathbb{R}^{N}\backslash\{0\}.
\end{aligned}
\end{equation*}
\end{theorem}
\begin{remark}
In \cite{Su2018FractionalLaplacian},
the authors set an open problem.
Our problem $(\mathcal{P})$ is a variant of the open problem.
\end{remark}

\begin{remark}
In order to study problem
$(\mathcal{P})$,
we must study problem
$(\mathcal{P}_{1})$
as follows:
$$
-\Delta u
=
\sum_{i=1}^{k}
\left(
\int_{\mathbb{R}^{N}}
\frac{|u|^{2^{*}_{\alpha_{i}}}}{|x-y|^{\alpha_{i}}}
\mathrm{d}y
\right)
|u|^{2^{*}_{\alpha_{i}}-2}u
+
|u|^{2^{*}-2}u
,
\mathrm{~in~}
\mathbb{R}^{N},
\eqno(\mathcal{P}_{1})
$$
where the parameters are same to problem
$(\mathcal{P})$.
We need show the relation of critical value between problem
$(\mathcal{P})$
and
problem
$(\mathcal{P}_{1})$
as follows:
$$\tilde{c}_{0}>c_{0},$$
where
$\tilde{c}_{0}$
and
$c_{0}$ are defined in Section 2.
There are finite many Hardy--Littlewood--Sobolev critical upper exponents in
problem
$(\mathcal{P})$
and
problem
$(\mathcal{P}_{1})$,
it is difficult to show $\tilde{c}_{0}>c_{0}$.
By using Lemma \ref{lemma9}
and
$S=\tilde{S}
\left(
1-
\frac{4\zeta}{(N-2)^{2}}
\right)^{\frac{N-1}{N}}$,
we overcome this difficulty in Lemma \ref{lemma13}.

We point out that
$S=\tilde{S}
\left(
1-
\frac{4\zeta}{(N-2)^{2}}
\right)^{\frac{N-1}{N}}$
plays a key role in the proof of $\tilde{c}_{0}>c_{0}$.
However,
for fractional Laplacian operator,
this equation is unknown.
So, we could not apply our method to prove the open problem in \cite{Su2018FractionalLaplacian}.
\end{remark}

\begin{remark}
Problem
$(\mathcal{P})$ is invariant under the weighted dilation
$$u\mapsto \tau^{\frac{N-2}{2}}u(\tau x).$$
Therefore, the well known Mountain Pass theorem does not yield critical point,
but only the Palais--Smale sequence.
It is necessary to show the non--vanishing of Palais--Smale sequence.
There are finite many Hardy--Littlewood--Sobolev critical upper exponents in problem
$(\mathcal{P})$,
it is difficult to show the non--vanishing of Palais--Smale sequence.
By using Coulomb--Sobolev space, endpoint refined Sobolev inequality and Lemma \ref{lemma5},
we overcome this difficulty in Lemma \ref{lemma12}.
\end{remark}
This paper is organized as follows:
In Section 2, we present some notations.
In Section 3, we show some key lemmas.
In Section 4, we study the Nehari manifolds for
problem $(\mathcal{P})$
and
problem $(\mathcal{P}_{1})$.
In Section 5, we investigate the Palais--Smale sequence of Problem $(\mathcal{P})$.
In Section 6, we show $\tilde{c}_{0}>c_{0}$.
In Section 7, we show the proof of Theorem \ref{theorem1}.
\section{Preliminaries}
Recall that the space
$D^{1,2}(\mathbb{R}^{N})$
is the completion of
$C^{\infty}_{0}(\mathbb{R}^{N})$
with respect to the norm
$$
\|u\|_{D}^{2}=\int_{\mathbb{R}^{N}}|\nabla u|^{2}\mathrm{d}x.
$$
It is well known that
$\frac{(N-2)^{2}}{4}$
is the best constant in the Hardy inequality
$$
\frac{(N-2)^{2}}{4}
\int_{\mathbb{R}^{N}}
\frac{ u^{2}}{|x|^{2}}
\mathrm{d}x
\leqslant
\int_{\mathbb{R}^{N}}
|\nabla u|^{2}
\mathrm{d}x
,~~
\mathrm{for~any~}
u\in
D^{1,2}(\mathbb{R}^{N}).
$$
By Hardy inequality and $\zeta\in(0,\frac{(N-2)^{2}}{4})$,
we derive that
$$
\|u\|_{\zeta}^{2}
=
\int_{\mathbb{R}^{N}}
\left(
|\nabla u|^{2}
-
\zeta
\frac{ u^{2}}{|x|^{2}}
\right)
\mathrm{d}x,
$$
is an equivalent norm in
$D^{1,2}(\mathbb{R}^{N})$,
and the following inequalities hold:
$$\left(1-\frac{4\zeta}{(N-2)^{2}}\right)\|u\|_{D}^{2}\leqslant\|u\|_{\zeta}^{2}\leqslant\|u\|_{D}^{2}.$$
For
$\alpha\in(0,N)$,
the Coulomb--Sobolev space \cite{Mercuri2016} is defined by
$$
\mathcal{E}^{1,\alpha,2^{*}_{\alpha}}(\mathbb{R}^{N})=
\left\{
\|u\|_{D}<\infty
~\mathrm{and}~
\int_{\mathbb{R}^{N}}
\int_{\mathbb{R}^{N}}
\frac{|u(x)|^{2^{*}_{\alpha}}|u(y)|^{2^{*}_{\alpha}}}{|x-y|^{\alpha}}
\mathrm{d}x
\mathrm{d}y<\infty
\right\}.
$$
We endow the space
$\mathcal{E}^{1,\alpha,2^{*}_{\alpha}}(\mathbb{R}^{N})$
with the norm
$$
\|u\|_{\mathcal{E},\alpha}^{2}
=
\|u\|_{D}^{2}
+
\left(
\int_{\mathbb{R}^{N}}
\int_{\mathbb{R}^{N}}
\frac{|u(x)|^{2^{*}_{\alpha}}|u(y)|^{2^{*}_{\alpha}}}{|x-y|^{\alpha}}
\mathrm{d}x
\mathrm{d}y
\right)^{\frac{1}{2^{*}_{\alpha}}}.
$$
For
$\alpha\in(0,N)$
and $\zeta\in(0,\frac{(N-2)^{2}}{4})$,
we could define the best constants:
\begin{equation}\label{10}
\begin{aligned}
S:=
\inf_{u\in D^{1,2}(\mathbb{R}^{N})\setminus\{0\}}
\frac{\|u\|^{2}_{\zeta}}
{(\int_{\mathbb{R}^{N}}
|u|^{2^{*}}
\mathrm{d}x)^{\frac{2}{2^{*}}}},
\end{aligned}
\end{equation}
and
\begin{equation}\label{11}
\begin{aligned}
S_{\alpha}:=
\inf_{u\in D^{1,2}(\mathbb{R}^{N})\setminus\{0\}}
\frac{\|u\|^{2}_{\zeta}}
{\left(
\int_{\mathbb{R}^{N}}
\int_{\mathbb{R}^{N}}
\frac{|u(x)|^{2^{*}_{\alpha}}|u(y)|^{2^{*}_{\alpha}}}{|x-y|^{\alpha}}
\mathrm{d}x
\mathrm{d}y
\right)^{\frac{1}{2^{*}_{\alpha}}}},
\end{aligned}
\end{equation}
and
\begin{equation}\label{12}
\begin{aligned}
\tilde{S}:=
\inf_{u\in D^{1,2}(\mathbb{R}^{N})\setminus\{0\}}
\frac{\|u\|^{2}_{D}}
{(\int_{\mathbb{R}^{N}}
|u|^{2^{*}}
\mathrm{d}x)^{\frac{2}{2^{*}}}},
\end{aligned}
\end{equation}
and
\begin{equation}\label{13}
\begin{aligned}
\tilde{S}_{\alpha}:=
\inf_{u\in D^{1,2}(\mathbb{R}^{N})\setminus\{0\}}
\frac{\|u\|^{2}_{D}}
{\left(
\int_{\mathbb{R}^{N}}
\int_{\mathbb{R}^{N}}
\frac{|u(x)|^{2^{*}_{\alpha}}|u(y)|^{2^{*}_{\alpha}}}{|x-y|^{\alpha}}
\mathrm{d}x
\mathrm{d}y
\right)^{\frac{1}{2^{*}_{\alpha}}}},
\end{aligned}
\end{equation}
where
$S$
and
$S_{\alpha}$
are attained in
$\mathbb{R}^{N}$
(see \cite[Lemma 2.1]{Yang2017}),
and
$\tilde{S}$
and
$\tilde{S}_{\alpha}$
are attained in
$\mathbb{R}^{N}$
(see \cite[Lemma 1.2]{Gao2016}).
Furthermore,
we know (see \cite[Formula (1)]{Swanson1992})
\begin{equation}\label{14}
\begin{aligned}
\tilde{S}=
\pi N (N-2)
\left(
\frac{\Gamma(\frac{N}{2})}{\Gamma(N)}
\right)^{\frac{2}{N}}.
\end{aligned}
\end{equation}
And see \cite[Lemma 1.2]{Gao2016}
\begin{equation}\label{15}
\begin{aligned}
\tilde{S}_{\alpha}=
\frac{\tilde{S}}{C(N,\alpha)^{\frac{1}{2^{*}_{\alpha}}}}.
\end{aligned}
\end{equation}
And see \cite[Formula (7)]{Felli2006}
\begin{equation}\label{16}
\begin{aligned}
S=\tilde{S}
\left(
1-
\frac{4\zeta}{(N-2)^{2}}
\right)^{\frac{N-1}{N}}.
\end{aligned}
\end{equation}
A measurable function
$u:\mathbb{R}^{N}\rightarrow \mathbb{R}$
belongs to the Morrey space
$\|u\|_{\mathcal{L}^{p,\varpi}}(\mathbb{R}^{N})$
with
$p\in[1,\infty)$
and
$\varpi\in(0,N]$
if and only if
$$
\|u\|^{p}_{\mathcal{L}^{p,\varpi}(\mathbb{R}^{N})}
=
\sup_{R>0,x\in\mathbb{R}^{N}}
R^{\varpi-N}
\int_{B(x,R)}
|u(y)|^{p}
\mathrm{d}y
<\infty.
$$
\begin{lemma}
\label{lemma4}
$\left.\right.$
\cite{Palatucci2014}
Let
$N\geqslant3$.
There exists
$C_{2}>0$
such that
for
$\iota$
and
$\vartheta$
satisfying
$\frac{2}{2^{*}}\leqslant\iota<1$,
$1\leqslant \vartheta<2^{*}=\frac{2N}{N-2}$,
we have
\begin{align*}
\left(
\int_{\mathbb{R}^{N}}
|u|^{2^{*}}
\mathrm{d}x
\right)^{\frac{1}{2^{*}}}
\leqslant
C_{2}
\|u\|_{D}^{\iota}
\|u\|_{\mathcal{L}^{\vartheta,\frac{\vartheta(N-2)}{2}}(\mathbb{R}^{N})}^{1-\iota},
\end{align*}
for any
$u\in D^{1,2}(\mathbb{R}^{N})$.
\end{lemma}
We introduce the energy functional associated to problem $(\mathcal{P})$ by
\begin{equation*}
\begin{aligned}
I(u)
=&
\frac{1}{2}
\|u\|_{\zeta}^{2}
-
\sum_{i=1}^{k}
\frac{1}{2\cdot2^{*}_{\alpha_{i}}}
\int_{\mathbb{R}^{N}}
\int_{\mathbb{R}^{N}}
\frac{|u(x)|^{2^{*}_{\alpha_{i}}}|u(y)|^{2^{*}_{\alpha_{i}}}}{|x-y|^{\alpha_{i}}}
\mathrm{d}x
\mathrm{d}y
-
\frac{1}{2^{*}}
\int_{\mathbb{R}^{N}}
|u|^{2^{*}}
\mathrm{d}x.
\end{aligned}
\end{equation*}
We also introduce the energy functional associated to problem $(\mathcal{P}_{1})$ by
\begin{equation*}
\begin{aligned}
\tilde{I}_{0}(u)
=&
\frac{1}{2}
\|u\|_{D}^{2}
-
\sum_{i=1}^{k}
\frac{1}{2\cdot2^{*}_{\alpha_{i}}}
\int_{\mathbb{R}^{N}}
\int_{\mathbb{R}^{N}}
\frac{|u(x)|^{2^{*}_{\alpha_{i}}}|u(y)|^{2^{*}_{\alpha_{i}}}}{|x-y|^{\alpha_{i}}}
\mathrm{d}x
\mathrm{d}y
-
\frac{1}{2^{*}}
\int_{\mathbb{R}^{N}}
|u|^{2^{*}}
\mathrm{d}x.
\end{aligned}
\end{equation*}
The Nehari manifold associated with
problem $(\mathcal{P})$ is  defined by
$$\mathcal{N}=\{u\in D^{1,2}(\mathbb{R}^{N})|\langle I^{'}(u),u\rangle=0,~u\not=0 \},$$
and
$$c_{0}=\inf_{u\in\mathcal{N}}I(u)~,~c_{1}=\inf_{u\in D^{1,2}(\mathbb{R}^{N})}\max_{t\geqslant0}I(tu)~\mathrm{and}~c=\inf_{\Upsilon\in\Gamma}
\max_{t\in [0,1]}
I(\Upsilon(t)),$$
where
$
\Gamma=
\{
\Upsilon\in C([0,1],D^{1,2}(\mathbb{R}^{N}))
:
\Upsilon(0)=0,
I(\Upsilon(1))<0
\}
$.

The Nehari manifold associated with
problem $(\mathcal{P}_{1})$ is  defined by
$$\widetilde{\mathcal{N}}=\{u\in D^{1,2}(\mathbb{R}^{N})|\langle \tilde{I}_{0}^{'}(u),u\rangle=0,~u\not=0 \},$$
and
$$
\tilde{c}_{0}=\inf_{u\in\widetilde{\mathcal{N}}}\tilde{I}_{0}(u)
~,~
\tilde{c}_{1}=\inf_{u\in D^{1,2}(\mathbb{R}^{N})}\max_{t\geqslant0}\tilde{I}_{0}(tu)
~\mathrm{and}~
\tilde{c}=\inf_{\tilde{\Upsilon}\in\tilde{\Gamma}}
\max_{t\in [0,1]}
\tilde{I}_{0}(\tilde{\Upsilon}(t)),$$
where
$
\tilde{\Gamma}=
\{
\tilde{\Upsilon}\in C([0,1],D^{1,2}(\mathbb{R}^{N}))
:
\tilde{\Upsilon}(0)=0,
\tilde{I}_{0}(\tilde{\Upsilon}(1))<0
\}
$.

\section{Some key Lemmas}
We show some properties of Coulomb--Sobolev space
$\mathcal{E}^{1,\alpha,2^{*}_{\alpha}}(\mathbb{R}^{N})$.
\begin{lemma}\label{lemma5}
Let
$( H_{1})$
hold.
If
$u\in \mathcal{E}^{1,\alpha_{j},2^{*}_{\alpha_{j}}}(\mathbb{R}^{N})$
$(j=1,\ldots,k)$,
then

\noindent
(i) $\|\cdot\|_{D}$
is an equivalent norm in
$\mathcal{E}^{1,\alpha_{j},2^{*}_{\alpha_{j}}}(\mathbb{R}^{N})$;

\noindent
(ii)
$
u\in
\bigcap_{i=1,i\not=j}^{k}
\mathcal{E}^{1,\alpha_{i},2^{*}_{\alpha_{i}}}(\mathbb{R}^{N})$;

\noindent
(iii)
$\|\cdot\|_{\mathcal{E},\alpha_{i}}$
are equivalent norms in
$\mathcal{E}^{1,\alpha_{j},2^{*}_{\alpha_{j}}}(\mathbb{R}^{N})$,
where
$i\not=j$
and
$i=1,\ldots,k$.
\end{lemma}
\noindent
{\bf Proof.}
\noindent
{\bf (1).}
Set
$j=1,\ldots,k$.
For any $u\in\mathcal{E}^{1,\alpha_{j},2^{*}_{\alpha_{j}}}(\mathbb{R}^{N})$,
applying the definition of Coulomb--Sobolev space,
we know
\begin{align}\label{17a}
\|u\|_{D}^{2}
\leqslant
\|u\|_{\mathcal{E},\alpha_{j}}^{2}
<\infty.
\end{align}
This implies that
$\mathcal{E}^{1,\alpha_{j},2^{*}_{\alpha_{j}}}(\mathbb{R}^{N})\subset D^{1,2}(\mathbb{R}^{N})$.
According to  $\mathcal{E}^{1,\alpha_{j},2^{*}_{\alpha_{j}}}(\mathbb{R}^{N})\subset D^{1,2}(\mathbb{R}^{N})$
and \eqref{13},
we have
\begin{align}\label{17b}
\|u\|_{\mathcal{E},\alpha_{j}}^{2}
\leqslant
\left(
1
+
\frac{1}{\tilde{S}_{\alpha_{j}}}
\right)
\|u\|_{D}^{2}.
\end{align}
Combining
\eqref{17a}
and
\eqref{17b},
we obtain
\begin{align}\label{17c}
\|u\|_{D}^{2}
\leqslant
\|u\|_{\mathcal{E},\alpha_{j}}^{2}
\leqslant
\left(
1
+
\frac{1}{\tilde{S}_{\alpha_{j}}}
\right)
\|u\|_{D}^{2}.
\end{align}
These imply that
$\|\cdot\|_{D}$
is an equivalent norm in
$\mathcal{E}^{1,\alpha_{j},2^{*}_{\alpha_{j}}}(\mathbb{R}^{N})$.

\noindent
{\bf (2).}
For any $u\in\mathcal{E}^{1,\alpha_{j},2^{*}_{\alpha_{j}}}(\mathbb{R}^{N})\subset D^{1,2}(\mathbb{R}^{N})$,
by using \eqref{17a}
and
\eqref{13},
we know
\begin{align}\label{17d}
\tilde{S}_{\alpha_{i}}
\left(
\int_{\mathbb{R}^{N}}
\int_{\mathbb{R}^{N}}
\frac{|u(x)|^{2^{*}_{\alpha_{i}}}|u(y)|^{2^{*}_{\alpha_{i}}}}{|x-y|^{\alpha_{i}}}
\mathrm{d}x
\mathrm{d}y
\right)^{\frac{1}{2^{*}_{\alpha_{i}}}}
\leqslant
\|u\|_{D}^{2}
\leqslant
\|u\|_{\mathcal{E},\alpha_{j}}^{2}
<\infty,
\end{align}
where
$i\not=j$
and
$i=1,\ldots,k$.
The inequality
\eqref{17d}
gives that
\begin{align*}
\|u\|_{\mathcal{E},\alpha_{i}}^{2}
=
\|u\|_{D}^{2}
+
\left(
\int_{\mathbb{R}^{N}}
\int_{\mathbb{R}^{N}}
\frac{|u(x)|^{2^{*}_{\alpha_{i}}}|u(y)|^{2^{*}_{\alpha_{i}}}}{|x-y|^{\alpha_{i}}}
\mathrm{d}x
\mathrm{d}y
\right)^{\frac{1}{2^{*}_{\alpha_{i}}}}
<\infty.
\end{align*}
This implies that $u\in\bigcap_{i=1,i\not=j}^{k}\mathcal{E}^{1,\alpha_{i},2^{*}_{\alpha_{i}}}(\mathbb{R}^{N})$.

\noindent
{\bf (3).}
For any $u\in\mathcal{E}^{1,\alpha_{j},2^{*}_{\alpha_{j}}}(\mathbb{R}^{N})$,
by using \eqref{17b},
we have
\begin{equation*}
\begin{aligned}
\|u\|_{\mathcal{E},\alpha_{j}}^{2}
\leqslant
\left(
\frac{\tilde{S}_{\alpha_{j}}+1}{\tilde{S}_{\alpha_{j}}}
\right)
\|u\|_{D}^{2}
\leqslant
\left(
\frac{\tilde{S}_{\alpha_{j}}+1}{\tilde{S}_{\alpha_{j}}}
\right)
\|u\|_{\mathcal{E},\alpha_{i}}^{2},
\end{aligned}
\end{equation*}
which imply that
\begin{equation*}
\begin{aligned}
\left(
\frac{\tilde{S}_{\alpha_{j}}}{\tilde{S}_{\alpha_{j}}+1}
\right)
\|u\|_{\mathcal{E},\alpha_{i}}^{2}
\leqslant
\|u\|_{\mathcal{E},\alpha_{j}}^{2}
\leqslant
\left(
\frac{\tilde{S}_{\alpha_{j}}+1}{\tilde{S}_{\alpha_{j}}}
\right)
\|u\|_{\mathcal{E},\alpha_{i}}^{2},
\end{aligned}
\end{equation*}
where
$0<\frac{\tilde{S}_{\alpha_{j}}}{S_{\alpha_{j}}+1}
<1<
\frac{\tilde{S}_{\alpha_{j}}+1}{\tilde{S}_{\alpha_{j}}}<\infty
$.
\qed

The following result is  the refinement of Hardy-Littlewood-Sobolev inequality.
\begin{lemma}\label{lemma6}
For any
$\alpha\in(0,N)$,
there exists
$C_{3}>0$
such that
for
$\iota$
and
$\vartheta$
satisfying
$\frac{2}{2^{*}}\leqslant\iota<1$,
$1\leqslant \vartheta<2^{*}=\frac{2N}{N-2}$,
we have
\begin{align*}
\left(
\int_{\mathbb{R}^{N}}
\int_{\mathbb{R}^{N}}
\frac{|u(x)|^{2^{*}_{\alpha}}|u(y)|^{2^{*}_{\alpha}}}{|x-y|^{\alpha}}
\mathrm{d}x
\mathrm{d}y
\right)^{\frac{1}{2^{*}_{\alpha}}}
\leqslant
C_{3}
\|u\|_{D}^{2\iota}
\|u\|_{\mathcal{L}^{\vartheta,\frac{\vartheta(N-2)}{2}}(\mathbb{R}^{N})}^{2(1-\iota)},
\end{align*}
for any
$u\in D^{1,2}(\mathbb{R}^{N})$.
\end{lemma}
\noindent
{\bf Proof.}
Let
$\frac{2}{2^{*}}\leqslant\iota<1$
and
$1\leqslant \vartheta<2^{*}=\frac{2N}{N-2}$.
By Hardy-Littlewood-Sobolev inequality
and
Lemma \ref{lemma4},
we have
\begin{equation*}
\begin{aligned}
\left(
\int_{\mathbb{R}^{N}}
\int_{\mathbb{R}^{N}}
\frac{|u(x)|^{2^{*}_{\alpha}}|u(y)|^{2^{*}_{\alpha}}}{|x-y|^{\alpha}}
\mathrm{d}x
\mathrm{d}y
\right)^{\frac{1}{2^{*}_{\alpha}}}
\leqslant&
C(N,\alpha)^{\frac{1}{2^{*}_{\alpha}}}
\|u\|_{L^{2^{*}}(\mathbb{R}^{N})}^{2}\\
\leqslant&
C(N,\alpha)^{\frac{1}{2^{*}_{\alpha}}}
C_{2}^{2}
\|u\|_{D}^{2\iota}
\|u\|_{\mathcal{L}^{\vartheta,\frac{\vartheta(N-2)}{2}}(\mathbb{R}^{N})}^{2(1-\iota)}.
\end{aligned}
\end{equation*}
\qed
\section{Nehari manifolds for problem $(\mathcal{P})$ and problem $(\mathcal{P}_{1})$}
We prove some properties of the Nehari manifold associated with
problem $(\mathcal{P})$.
\begin{lemma}\label{lemma7}
Let
$N\geqslant3$,
$\zeta\in(0,\frac{(N-2)^{2}}{4})$
and
$(H_{1})$
hold.
Then
$$c_{0}=\inf_{u\in\mathcal{N}}I(u)>0.$$
\end{lemma}
\noindent
{\bf Proof.}
We divide our proof into two steps.

\noindent
{\bf Step 1.}
We claim that
any limit point of a sequence in
$\mathcal{N}$
is different from zero.
According to
$\langle I^{'}(u),u\rangle=0$,
(\ref{10})
and
(\ref{11}),
for any $u\in\mathcal{N}$,
we obtain
\begin{equation*}
\begin{aligned}
0
=
\langle I^{'}(u),u\rangle
\geqslant&
\|u\|^{2}_{\zeta}
-
\frac{1}{S ^{\frac{2^{*} }{2}}}
\|u\|_{\zeta}^{2^{*} }
-
\sum_{i=1}^{k}
\frac{1}{S_{\alpha_{i}} ^{2^{*}_{\alpha_{i}}}}
\|u\|_{\zeta}^{2\cdot 2^{*}_{\alpha_{i}} }.
\end{aligned}
\end{equation*}
From above expression,
we have
\begin{equation}\label{19}
\begin{aligned}
\|u\|^{2}_{\zeta}
\leqslant
\frac{1}{S ^{\frac{2^{*} }{2}}}
\|u\|_{\zeta}^{2^{*} }
+
\sum_{i=1}^{k}
\frac{1}{S_{\alpha_{i}} ^{2^{*}_{\alpha_{i}}}}
\|u\|_{\zeta}^{2\cdot 2^{*}_{\alpha_{i}} }.
\end{aligned}
\end{equation}
Set
$$\kappa
:=
\frac{1}{S ^{\frac{2^{*} }{2}}}
+
\sum_{i=1}^{k}
\frac{1}
{S_{\alpha_{i}} ^{2^{*}_{\alpha_{i}}}}
.$$
Applying
(\ref{10}),
(\ref{11})
and
$(H_{1})$,
we get
$$
0<\kappa<\infty.$$
From
$(H_{1})$,
we know
$$\frac{2N-\alpha_{1}}{N-2}=\frac{2N-\alpha_{i}}{N-2}.$$
Now the proof of Step 1 is divided into two cases:
(i)
$\|u\|_{\zeta}\geqslant1$;
(ii)
$\|u\|_{\zeta}<1$.

\noindent
{\bf Case (i).}
$\|u\|_{\zeta}\geqslant1$.
The inequality
(\ref{19})
gives
$$
\|u\|_{\zeta}^{2}
\leqslant
\frac{1}{S ^{\frac{2^{*} }{2}}}
\|u\|_{\zeta}^{2^{*} }
+
\sum_{i=1}^{k}
\frac{1}{S_{\alpha_{i}} ^{2^{*}_{\alpha_{i}}}}
\|u\|_{\zeta}^{2\cdot 2^{*}_{\alpha_{i}} }
\leqslant
\kappa
\|u\|_{\zeta}^{2\cdot2^{*}_{\alpha_{1}}},
$$
which implies that
\begin{equation}\label{20}
\begin{aligned}
\|u\|_{\zeta}
\geqslant
\kappa
^{\frac{1}{2-2\cdot2^{*}_{\alpha_{1}}}}
.
\end{aligned}
\end{equation}

\noindent
{\bf Case (ii).}
$\|u\|_{\zeta}<1$.
Again,
by
(\ref{19}),
we know
\begin{equation}\label{21}
\begin{aligned}
\|u\|_{\zeta}
\geqslant
\kappa
^{\frac{1}{2-2^{*} }}
.
\end{aligned}
\end{equation}
According to
(\ref{20})
and
(\ref{21}),
we deduce that
\begin{equation}\label{22}
\begin{aligned}
\|u\|_{\zeta}
\geqslant
\min\left\{
\kappa
^{\frac{1}{2-2\cdot2^{*}_{\alpha_{1}}}}
,
\kappa
^{\frac{1}{2-2^{*} }}
\right\}.
\end{aligned}
\end{equation}
Hence,
we know that
any limit point of a sequence in $\mathcal{N}$ is different from zero.

\noindent
{\bf Step 2.}
Now,
we claim that
$I$
is bounded from below on
$\mathcal{N}$.
For any
$u\in \mathcal{N}$,
by using
(\ref{22}),
we get
\begin{equation*}
\begin{aligned}
I(u)
\geqslant
\left(
\frac{1}{2}
-
\frac{1}{2^{*}}
\right)
\|u\|^{2}_{\zeta}
\geqslant
\frac{1}{N}
\min\left\{
\kappa
^{\frac{2}{2-2\cdot2^{*}_{\alpha_{1}}}}
,
\kappa
^{\frac{2}{2-2^{*} }}
\right\}.
\end{aligned}
\end{equation*}
Therefore,
$I$
is bounded from below on
$\mathcal{N}$,
and
$c_{0}>0$.
\qed
\begin{lemma}\label{lemma8}
Let
$N\geqslant3$,
$\zeta\in(0,\frac{(N-2)^{2}}{4})$
and
$(H_{1})$
hold.
Then

\noindent
(i)
for each
$u\in D^{1,2}(\mathbb{R}^{N})\setminus\{0\}$,
there exists a unique
$t_{u}>0$
such that
$t_{u}u\in \mathcal{N}$;

\noindent
(ii)
$c_{0}=c_{1}=c>0$.
\end{lemma}
\noindent
{\bf Proof.}
The proof is standard, so we sketch it. Further details can be derived as in the proofs of Theorem 4.1 and 4.2 in
\cite{Willem1996}.
We omit it.
\qed

We prove some properties of the Nehari manifold associated with
problem $(\mathcal{P}_{1})$.
\begin{lemma}\label{lemma9}
Let
$N\geqslant3$
and
$(H_{1})$
hold.
For any
$u\in \widetilde{\mathcal{N}}$,
we have
\begin{equation*}
\begin{aligned}
\|u\|_{D}
\geqslant
\min\left\{
\tilde{\kappa}^{\frac{1}{2-2\cdot2^{*}_{\alpha_{1}}}}
,
\tilde{\kappa}
^{\frac{1}{2-2^{*} }}
\right\},
\end{aligned}
\end{equation*}
where
$$\tilde{\kappa}
=
\frac{1}{\tilde{S} ^{\frac{2^{*} }{2}}}
+
\sum_{i=1}^{k}
\frac{1}
{\tilde{S}_{\alpha_{i}} ^{2^{*}_{\alpha_{i}}}}
.$$
And,
$$\tilde{c}_{1}=\tilde{c}=\tilde{c}_{0}=\inf_{u\in\widetilde{\mathcal{N}}}\tilde{I}_{0}(u)>0.$$
\end{lemma}
\section{Analysis of the Palais--Smale sequence for Problem $(\mathcal{P})$}
We show that the functional $I$ satisfies the Mountain Pass geometry, and
estimate the Mountain Pass level.
\begin{lemma}\label{lemma11}
Let
$N\geqslant3$,
$\zeta\in(0,\frac{(N-2)^{2}}{4})$
and
$(H_{1})$
hold.
Then
there exists a
$(PS)_{c}$
sequence of
$I$
at level
$c$,
where
$$0<c<c^{*}=
\min
\left\{
\frac
{N+2-\alpha_{1}}
{2(2N-\alpha_{1})}
S_{\alpha_{1}}^{\frac{2N-\alpha_{1}}{N+2-\alpha_{1}}}
,
\ldots
,
\frac{N+2-\alpha_{k}}
{2(2N-\alpha_{k})}
S_{\alpha_{k}}^{\frac{2N-\alpha_{k}}{N+2-\alpha_{k}}}
,
\frac{1}{N}S ^{\frac{N}{2}}
\right\}.$$
\end{lemma}
\noindent
{\bf Proof.}
We divide our proof into two steps.

\noindent
{\bf Step 1.}
We prove that
$I$
satisfies all the conditions in Mountain Pass theorem.

\noindent
(i)
$I(0)=0$;

\noindent
(ii)
For any
$u\in D^{1,2}(\mathbb{R}^{N})\setminus\{0\}$,
we have
\begin{equation*}
\begin{aligned}
I(u)
\geqslant
\frac{1}{2}
\|u\|_{\zeta}^{2}
-
\sum_{i=1}^{k}
\frac{1}{2\cdot2^{*}_{\alpha_{i}}S_{\alpha_{1}}^{2^{*}_{\alpha_{i}}}}
\|u\|_{\zeta}^{2\cdot 2^{*}_{\alpha_{i}}}
-
\frac{1}{2^{*}S^{\frac{2^{*}}{2}}}
\|u\|_{\zeta}^{2^{*}}
.
\end{aligned}
\end{equation*}
Because of
$2<2^{*}<2\cdot 2^{*}_{\alpha_{k}}<\cdots<2\cdot 2^{*}_{\alpha_{1}}$,
there exists a sufficiently small positive number
$\rho$
such that
$$\varsigma:=\inf_{\|u\|_{\zeta}=\rho}I(u)>0=I(0).$$

\noindent
(iii)
Given
$u\in D^{1,2}(\mathbb{R}^{N})\setminus\{0\}$
such that
$\lim_{t\rightarrow\infty}
I(tu)=-\infty.$
We choose
$t_{u}>0$
corresponding to
$u$
such that
$I(tu)<0$
for all
$t>t_{u}$
and
$\|t_{u} u\|_{\zeta}>\rho$.
Set
$$c=
\inf_{\Upsilon\in\Gamma}
\max_{t\in [0,1]}
I(\Upsilon(t)),$$
where
$
\Gamma=
\{
\Upsilon\in C([0,1],D^{1,2}(\mathbb{R}^{N}))
:
\Upsilon(0)=0,
\Upsilon(1)=t_{u}u
\}
$.

\noindent
{\bf Step 2.}
Here we show
$0<c<c^{*}
$.
By using
(\ref{10})
and
(\ref{11}),
there exist minimizers
$u_{i}\in D^{1,2}(\mathbb{R}^{N})\setminus\{0\}$
of
$S_{\alpha_{i}}$
and
$\tilde{u}\in D^{1,2}(\mathbb{R}^{N})\setminus\{0\}$
of
$S$.
For
$t\geqslant 0$,
we set
\begin{equation*}
\begin{aligned}
f_{i}(t)
=&
\frac{t^{2}}{2}
\|u_{i}\|_{\zeta}^{2}
-
\frac{ t^{2\cdot2^{*}_{\alpha_{i}}}}{2\cdot2^{*}_{\alpha_{i}}}
\int_{\mathbb{R}^{N}}
\int_{\mathbb{R}^{N}}
\frac{|u_{i}(x)|^{2^{*}_{\alpha_{i}}}|u_{i}(y)|^{2^{*}_{\alpha_{i}}}}{|x-y|^{\alpha_{i}}}
\mathrm{d}x
\mathrm{d}y,
~\mathrm{for}~
i=1,\ldots,k,
\end{aligned}
\end{equation*}
and
\begin{equation*}
\begin{aligned}
\tilde{f}(t)
=&
\frac{t^{2}}{2}
\|\tilde{u}\|_{\zeta}^{2}
-
\frac{t^{2^{*}}}{2^{*}}
\int_{\mathbb{R}^{N}}
|\tilde{u}|^{2^{*}}
\mathrm{d}x
.
\end{aligned}
\end{equation*}
If
$
\frac{N+2-\alpha_{i}}
{2(2N-\alpha_{i})}
S_{\alpha_{i}}^{\frac{2N-\alpha_{i}}{N+2-\alpha_{i}}}
=
\min
\left\{
\frac{N+2-\alpha_{1}}
{2(2N-\alpha_{1})}
S_{\alpha_{1}}^{\frac{2N-\alpha_{1}}{N+2-\alpha_{1}}}
,
\ldots
,
\frac{N+2-\alpha_{k}}
{2(2N-\alpha_{k})}
S_{\alpha_{k}}^{\frac{2N-\alpha_{k}}{N+2-\alpha_{k}}}
,
\frac{1}{N}S ^{\frac{N}{2}}
\right\},
$
for
$i=1,\ldots,k$,
we know that
$f_{i}^{'}(\cdot)=0$
if and only if
\begin{equation*}
\begin{aligned}
t
\|u_{i}\|_{\zeta}^{2}
=
t^{2\cdot2^{*}_{\alpha_{i}}-1}
\int_{\mathbb{R}^{N}}
\int_{\mathbb{R}^{N}}
\frac{|u_{i}(x)|^{2^{*}_{\alpha_{i}}}|u_{i}(y)|^{2^{*}_{\alpha_{i}}}}{|x-y|^{\alpha_{i}}}
\mathrm{d}x
\mathrm{d}y,
~\mathrm{for}~
i=1,\ldots,k.
\end{aligned}
\end{equation*}
We can see that
$f_{i}(\cdot)$ achieve their maximums at
$t_{i}$
as follows:
\begin{equation*}
\begin{aligned}
t_{i}^{2\cdot2^{*}_{\alpha_{i}}-2}
=
\frac{\|u_{i}\|_{\zeta}^{2}}
{\int_{\mathbb{R}^{N}}
\int_{\mathbb{R}^{N}}
\frac{|u_{i}(x)|^{2^{*}_{\alpha_{i}}}|u_{i}(y)|^{2^{*}_{\alpha_{i}}}}{|x-y|^{\alpha_{i}}}
\mathrm{d}x
\mathrm{d}y}
,
~\mathrm{for}~
i=1,\ldots,k.
\end{aligned}
\end{equation*}
Thus we obtain
\begin{equation*}
\begin{aligned}
0<c\leqslant&
\sup_{t\geqslant0}
I(tu_{i})
\leqslant
\max_{t\geqslant0}
f_{i}(t)
=
f_{i}(t_{i})
=
\frac{N+2-\alpha_{i}}
{2(2N-\alpha_{i})}
S_{\alpha_{i}}^{\frac{2N-\alpha_{i}}{N+2-\alpha_{i}}}.
\end{aligned}
\end{equation*}
Similar to \cite{Pucci2009,Yang2017,ZhangBL2017},
we get
\begin{equation*}
\begin{aligned}
\sup_{t\geqslant0}
I(tu_{i})
<
\max_{t\geqslant0}
f_{i}(t),
\end{aligned}
\end{equation*}
which means
$0<c<
\frac{N+2-\alpha_{i}}
{2(2N-\alpha_{i})}
S_{\alpha_{i}}^{\frac{2N-\alpha_{i}}{N+2-\alpha_{i}}}$.

If
$
\frac{1}{N}S ^{\frac{N}{2}}
=
\min
\left\{
\frac{N+2-\alpha_{1}}
{2(2N-\alpha_{1})}
S_{\alpha_{1}}^{\frac{2N-\alpha_{1}}{N+2-\alpha_{1}}}
,
\ldots
,
\frac{N+2-\alpha_{k}}
{2(2N-\alpha_{k})}
S_{\alpha_{k}}^{\frac{2N-\alpha_{k}}{N+2-\alpha_{k}}}
,
\frac{1}{N}S ^{\frac{N}{2}}
\right\}
$,
similarly, we have
\begin{equation*}
\begin{aligned}
0<c\leqslant
\sup_{t\geqslant0}
I(t\tilde{u})
<
\tilde{f}(t)
\leqslant
\max_{t\geqslant0}
\tilde{f}(t)
=
\frac{1}{N}
S^{\frac{N}{2}},
\end{aligned}
\end{equation*}
which means
$0<c<\frac{1}{N}
S^{\frac{N}{2}}$.
\qed

The following result implies the non--vanishing of
$(PS)_{c}$
sequence.
\begin{lemma}\label{lemma12}
Let
$N\geqslant3$,
$\zeta\in(0,\frac{(N-2)^{2}}{4})$
and
$(H_{1})$
hold.
Let
$\{u_{n}\}$
be a
$(PS)_{c}$
sequence of
$I$
with
$c\in(0,c^{*})$.
Then
\begin{equation*}
\begin{aligned}
\lim_{n\rightarrow\infty}
\int_{\mathbb{R}^{N}}
|u_{n}|^{2^{*}}
\mathrm{d}x>0,
\end{aligned}
\end{equation*}
and
\begin{equation*}
\begin{aligned}
\lim_{n\rightarrow\infty}
\int_{\mathbb{R}^{N}}
\int_{\mathbb{R}^{N}}
\frac{|u_{n}(x)|^{2^{*}_{\alpha_{i}}}|u_{n}(y)|^{2^{*}_{\alpha_{i}}}}{|x-y|^{\alpha_{i}}}
\mathrm{d}x
\mathrm{d}y>0,
~(i=1,\ldots,k).
\end{aligned}
\end{equation*}
\end{lemma}
\noindent
{\bf Proof.}
It is easy to see that
$\{u_{n}\}$
is uniformly bounded in
$D^{1,2}(\mathbb{R}^{N})$.
We divide our proof into three cases:

\noindent
(1)
$
\lim_{n\rightarrow\infty}
\int_{\mathbb{R}^{N}}
\int_{\mathbb{R}^{N}}
\frac{|u_{n}(x)|^{2^{*}_{\alpha_{1}}}|u_{n}(y)|^{2^{*}_{\alpha_{1}}}}{|x-y|^{\alpha_{1}}}
\mathrm{d}x
\mathrm{d}y>0
$;

\noindent
(2)
$
\lim_{n\rightarrow\infty}
\int_{\mathbb{R}^{N}}
\int_{\mathbb{R}^{N}}
\frac{|u_{n}(x)|^{2^{*}_{\alpha_{i}}}|u_{n}(y)|^{2^{*}_{\alpha_{i}}}}{|x-y|^{\alpha_{i}}}
\mathrm{d}x
\mathrm{d}y>0
$,
$i=2,\ldots,k$;

\noindent
(3)
$
\lim_{n\rightarrow\infty}
\int_{\mathbb{R}^{N}}
|u_{n}|^{2^{*}}
\mathrm{d}x>0
$.

\noindent
{\bf Case 1.}
Suppose on the contrary that
\begin{equation}\label{23}
\begin{aligned}
\lim_{n\rightarrow\infty}
\int_{\mathbb{R}^{N}}
\int_{\mathbb{R}^{N}}
\frac{|u_{n}(x)|^{2^{*}_{\alpha_{1}}}|u_{n}(y)|^{2^{*}_{\alpha_{1}}}}{|x-y|^{\alpha_{1}}}
\mathrm{d}x
\mathrm{d}y=0.
\end{aligned}
\end{equation}
Since
$\{u_{n}\}$
is uniformly bounded in
$D^{1,2}(\mathbb{R}^{N})$,
there exists a constant
$0<C<\infty$
such that
$\|u_{n}\|_{D}\leqslant C$.
According to (\ref{13})
and the definition of
Coulomb--Sobolev space,
we obtain
$u_{n}\in \mathcal{E}^{1,\alpha_{1},2^{*}_{\alpha_{1}}}(\mathbb{R}^{N})$.
Applying Lemma \ref{lemma2},
we have
\begin{equation}\label{24}
\begin{aligned}
&
\lim_{n\rightarrow\infty}
\|u_{n}\|_{L^{2^{*}}(\mathbb{R}^{N})}\\
\leqslant&
C
\left(
\lim_{n\rightarrow\infty}
\int_{\mathbb{R}^{N}}
\int_{\mathbb{R}^{N}}
\frac{|u_{n}(x)|^{2^{*}_{\alpha_{1}}}|u_{n}(y)|^{2^{*}_{\alpha_{1}}}}{|x-y|^{\alpha_{1}}}
\mathrm{d}x
\mathrm{d}y
\right)
^{\frac{N-2}{N(N+2-\alpha_{1})}}
=0.
\end{aligned}
\end{equation}
Combining Hardy--Littlewood--Sobolev inequality
and
(\ref{24}),
for all
$i=2,\ldots,k$,
we know
\begin{equation}\label{25}
\begin{aligned}
\lim_{n\rightarrow\infty}
\int_{\mathbb{R}^{N}}
\int_{\mathbb{R}^{N}}
\frac{|u_{n}(x)|^{2^{*}_{\alpha_{i}}}|u_{n}(y)|^{2^{*}_{\alpha_{i}}}}{|x-y|^{\alpha_{i}}}
\mathrm{d}x
\mathrm{d}y
\leqslant&
C
\lim_{n\rightarrow\infty}
\|u_{n}\|_{L^{2^{*}}(\mathbb{R}^{N})}^{2\cdot2^{*}_{\alpha_{i}}}
=0.
\end{aligned}
\end{equation}
Owing to
(\ref{23})
--
(\ref{25})
and the definition of
$(PS)_{c}$
sequence,
we obtain
$$
c
+
o(1)
=
\frac{1}{2}
\|u_{n}\|_{\zeta}^{2},$$
and
$$o(1)
=
\|u_{n}\|_{\zeta}^{2}.
$$
These imply that
$c=0$,
which contradicts as $0<c$.

\noindent
{\bf Case 2.}
From Case 1,
we have
$u_{n}\in \mathcal{E}^{1,\alpha_{1},2^{*}_{\alpha_{1}}}(\mathbb{R}^{N})$.
Applying the result of (ii) in Lemma \ref{lemma5},
we know that
$u_{n}\in \bigcap_{i=2}^{k}
\mathcal{E}^{1,\alpha_{i},2^{*}_{\alpha_{i}}}(\mathbb{R}^{N})$.
Then we could use the
endpoint refined Sobolev inequality for parameters $\alpha_{i}$
$(i=2,\ldots,k)$.
Similar to Case 1,
we prove that
$$\lim_{n\rightarrow\infty}
\int_{\mathbb{R}^{N}}
\int_{\mathbb{R}^{N}}
\frac{|u_{n}(x)|^{2^{*}_{\alpha_{i}}}|u_{n}(y)|^{2^{*}_{\alpha_{i}}}}{|x-y|^{\alpha_{i}}}
\mathrm{d}x
\mathrm{d}y>0,
~(i=2,\ldots,k).$$

\noindent
{\bf Case 3.}
Suppose that
\begin{equation}\label{26}
\begin{aligned}
\lim_{n\rightarrow\infty}
\int_{\mathbb{R}^{N}}
|u_{n}|^{2^{*}}
\mathrm{d}x
=0.
\end{aligned}
\end{equation}
By using Lemma \ref{lemma1},
for all
$i=1,\ldots,k$,
we have
\begin{equation}\label{27}
\begin{aligned}
\lim_{n\rightarrow\infty}
\int_{\mathbb{R}^{N}}
\int_{\mathbb{R}^{N}}
\frac{|u_{n}(x)|^{2^{*}_{\alpha_{i}}}|u_{n}(y)|^{2^{*}_{\alpha_{i}}}}{|x-y|^{\alpha_{i}}}
\mathrm{d}x
\mathrm{d}y
\leqslant&
C
\lim_{n\rightarrow\infty}
\|u_{n}\|_{L^{2^{*}}(\mathbb{R}^{N})}^{2\cdot2^{*}_{\alpha_{i}}}=0,\\
\end{aligned}
\end{equation}
Applying
(\ref{26})
and
(\ref{27}),
we get
$$
c
+
o(1)
=
\frac{1}{2}
\|u_{n}\|_{\zeta}^{2},
$$
and
$$o(1)
=
\|u_{n}\|_{\zeta}^{2}.
$$
These imply that
$c=0$,
which contradicts as $0<c$.
\qed
\section{The proof of $\tilde{c}_{0}>c_{0}$}
In this section,
we show that
$\tilde{c}_{0}>c_{0}$.
\begin{lemma}\label{lemma13}
Let
$N\geqslant3$,
$(H_{1})$
and
$(H_{2})$
hold.
If
$$\zeta
\in
\left(
\frac{(N-2)^{2}}{4}
\left(
1-
\frac{1}
{
(k+1)^{\frac{N-2}{N-1}}
\tilde{S}^{\frac{N}{N-1}}
}
\right)
,
\frac{(N-2)^{2}}{4}
\right),$$
then $\tilde{c}_{0}>c_{0}$.
\end{lemma}

\noindent
{\bf Proof.}
We divide our proof into three steps.

\noindent
{\bf Step 1.}
In this step,
we show the property of
$\tilde{c}_{0}$.
From Lemma \ref{lemma9},
there exists
$\bar{\bar{u}}$
such that
$\tilde{I}_{0}(\bar{\bar{u}})=\tilde{c}_{0}$
and
$\bar{\bar{u}}\in \widetilde{\mathcal{N}}$.
Then
\begin{equation}\label{28}
\begin{aligned}
\tilde{c}_{0}
=&
\tilde{I}_{0}(\bar{\bar{u}})
-
\frac{1}{2^{*}}
\langle \tilde{I}_{0}^{'}(\bar{\bar{u}}),\bar{\bar{u}}\rangle\\
=&
\left(
\frac{1}{2}
-
\frac{1}{2^{*}}
\right)
\|\bar{\bar{u}}\|^{2}_{D}
+
\sum_{i=1}^{k}
\left(
\frac{1}{2^{*}}
-
\frac{1}{2\cdot 2^{*}_{\alpha_{i}}}
\right)
\int_{\mathbb{R}^{N}}
\int_{\mathbb{R}^{N}}
\frac{|\bar{\bar{u}}(x)|^{2^{*}_{\alpha_{i}}}|\bar{\bar{u}}(y)|^{2^{*}_{\alpha_{i}}}}{|x-y|^{\alpha_{i}}}
\mathrm{d}x
\mathrm{d}y\\
\geqslant&
\left(
\frac{1}{2}
-
\frac{1}{2^{*}}
\right)
\|\bar{\bar{u}}\|^{2}_{D}
\\
\geqslant&
\frac{1}{N}
\min
\left\{
\tilde{\kappa}
^{-\frac{N-2}{N+2-\alpha_{1}}}
,
\tilde{\kappa}
^{-\frac{N-2}{2}}
\right\}
.
\end{aligned}
\end{equation}
\noindent
{\bf Step 2.}
In this step,
we show some basic results.
Firstly,
according to (\ref{14}) and $N\geqslant3$,
we get
$$\tilde{S}>1.$$
Secondly,
we will show the following results:
$$\zeta\in
\left(
\frac{
(N-2)^{2}}{4}
\left(
1-
\frac{1}
{
(k+1)^{\frac{N-2}{N-1}}
}
\right)
,
\frac{(N-2)^{2}}{4}
\right),$$
and
$$\zeta\in
\left(
\frac{
(N-2)^{2}}{4}
\left(
1-
\frac{1}
{
\tilde{S}^{\frac{N}{N-1}}
}
\right)
,
\frac{(N-2)^{2}}{4}
\right).$$
By using $\tilde{S}>1$ and $k+1>1$,
we have
$$
\frac{1}
{
(k+1)^{\frac{N-2}{N-1}}
\tilde{S}^{\frac{N}{N-1}}
}
<
\frac{1}
{
(k+1)^{\frac{N-2}{N-1}}
}
~\mathrm{and}~
\frac{1}
{
(k+1)^{\frac{N-2}{N-1}}
\tilde{S}^{\frac{N}{N-1}}
}
<
\frac{1}
{
\tilde{S}^{\frac{N}{N-1}}
}.
$$
These imply that
$$
1-
\frac{1}
{
(k+1)^{\frac{N-2}{N-1}}
\tilde{S}^{\frac{N}{N-1}}
}
>
1
-
\frac{1}
{
(k+1)^{\frac{N-2}{N-1}}
},$$
and
$$
1
-
\frac{1}
{
(k+1)^{\frac{N-2}{N-1}}
\tilde{S}^{\frac{N}{N-1}}
}
>
1
-
\frac{1}
{
\tilde{S}^{\frac{N}{N-1}}
}.
$$
Hence,
we deduce that
$$
\frac{4\zeta}{(N-2)^{2}}
\in
\left(
1-
\frac{1}
{
(k+1)^{\frac{N-2}{N-1}}
\tilde{S}^{\frac{N}{N-1}}
}
,
1
\right)
\subset
\left(
1-
\frac{1}
{
(k+1)^{\frac{N-2}{N-1}}
}
,
1
\right),$$
and
$$
\frac{4\zeta}{(N-2)^{2}}
\in
\left(
1-
\frac{1}
{
(k+1)^{\frac{N-2}{N-1}}
\tilde{S}^{\frac{N}{N-1}}
}
,
1
\right)
\subset
\left(
1-
\frac{1}
{
\tilde{S}^{\frac{N}{N-1}}
}
,
1
\right).$$

\noindent
{\bf Step 3.}
Here we show
$\tilde{c}_{0}>c_{0}$.
The proof of this step is divided into four cases:

\noindent
(1)
$\tilde{S}
=
\min
\left\{
\tilde{S}
,
\tilde{S}_{\alpha_{1}}
,\ldots,\tilde{S}_{\alpha_{k}}
\right\}$
and
$\tilde{\kappa}\geqslant1$;

\noindent
(2)
$\tilde{S}
=
\min
\left\{
\tilde{S}
,
\tilde{S}_{\alpha_{1}}
,\ldots,\tilde{S}_{\alpha_{k}}
\right\}$
and
$\tilde{\kappa}\leqslant1$;

\noindent
(3)
$1<\tilde{S}_{\alpha_{j}}
=
\min
\left\{
\tilde{S}
,
\tilde{S}_{\alpha_{1}}
,\ldots,\tilde{S}_{\alpha_{k}}
\right\}$
and
$\tilde{\kappa}\geqslant1$;

\noindent
(4)
$1<\tilde{S}_{\alpha_{j}}
=
\min
\left\{
\tilde{S}
,
\tilde{S}_{\alpha_{1}}
,\ldots,\tilde{S}_{\alpha_{k}}
\right\}$
and
$\tilde{\kappa}\leqslant1$.

\noindent
{\bf Case
(1).}
Since
$\frac{2^{*}}{2}<2^{*}_{\alpha_{k}}<\cdots<2^{*}_{\alpha_{1}}$
and
$1<\tilde{S}
=
\min
\left\{
\tilde{S}
,
\tilde{S}_{\alpha_{1}}
,\ldots,\tilde{S}_{\alpha_{k}}
\right\}$,
we have
\begin{equation}\label{29}
\begin{aligned}
\tilde{\kappa}
=
\frac{1}{\tilde{S} ^{\frac{2^{*} }{2}}}
+
\sum_{i=1}^{k}
\frac{1}
{\tilde{S}_{\alpha_{i}} ^{2^{*}_{\alpha_{i}} }}
\leqslant
\frac{k+1}{\tilde{S} ^{\frac{2^{*} }{2}}}.
\end{aligned}
\end{equation}
Combining
(\ref{28}),
(\ref{29}),
$\tilde{\kappa}\geqslant1$
and
$-\frac{N-2}{N+2-\alpha_{1}}>-\frac{N-2}{2}$,
we get
\begin{equation}\label{30}
\begin{aligned}
\tilde{c}_{0}
\geqslant
\frac{1}{N}
\min
\left\{
\tilde{\kappa}
^{-\frac{N-2}{N+2-\alpha_{1}}}
,
\tilde{\kappa}
^{-\frac{N-2}{2}}
\right\}
=
\frac{1}{N}
\tilde{\kappa}
^{-\frac{N-2}{2}}
\geqslant
\frac{\tilde{S} ^{\frac{N}{2}}}{N}
(k+1)^{-\frac{N-2}{2}}.
\end{aligned}
\end{equation}
By using
$\zeta
\in
\left(
\frac{
(N-2)^{2}}{4}\left(1-\frac{1}{(k+1)^{\frac{N-2}{N-1}}}\right)
,
\frac{(N-2)^{2}}{4}
\right)$,
we obtain
\begin{equation*}
\begin{aligned}
k+1
<&
\left(
1-\frac{4}{(N-2)^{2}}\zeta
\right)^{-\frac{N-1}{N-2}}\\
=&
\left(
\frac{\tilde{S}}{
\tilde{S}
\left(
1-\frac{4}{(N-2)^{2}}\zeta
\right)^{\frac{N-1}{N}}}
\right)^{\frac{N}{N-2}}\\
=&
\left(
\frac{\tilde{S}}{
S}
\right)^{\frac{N}{N-2}},
\end{aligned}
\end{equation*}
which implies that
\begin{equation}\label{31}
\begin{aligned}
(k+1)^{-\frac{N-2}{2}}
\tilde{S}^{\frac{N}{2}}
>
S^{\frac{N}{2}}.
\end{aligned}
\end{equation}
Putting
(\ref{31})
into
(\ref{30}),
we know
\begin{equation}\label{32}
\begin{aligned}
\tilde{c}_{0}
>
\frac{S ^{\frac{N}{2}}}{N}.
\end{aligned}
\end{equation}
According to
(\ref{32}),
Lemma \ref{lemma8}
and
Lemma \ref{lemma11},
we obtain
\begin{equation*}
\begin{aligned}
\tilde{c}_{0}
>
\frac{S ^{\frac{N}{2}}}{N}
>
c=c_{0}>0.
\end{aligned}
\end{equation*}

\noindent
{\bf Case (2).}
Combining
(\ref{28}),
(\ref{29}),
$\tilde{\kappa}\leqslant1$
and
$-\frac{N-2}{N+2-\alpha_{1}}>-\frac{N-2}{2}$,
we get
\begin{equation}\label{33}
\begin{aligned}
\tilde{c}_{0}
\geqslant&
\frac{1}{N}
\min
\left\{
\tilde{\kappa}
^{-\frac{N-2}{N+2-\alpha_{1}}}
,
\tilde{\kappa}
^{-\frac{N-2}{2}}
\right\}\\
=&
\frac{\tilde{\kappa}
^{-\frac{N-2}{N+2-\alpha_{1}}}}{N}
\geqslant
\frac{\tilde{S} ^{\frac{N}{N+2-\alpha_{1}}}}{N}
(k+1)
^{-\frac{N-2}{N+2-\alpha_{1}}}.
\end{aligned}
\end{equation}
Since
$\zeta\in
\left(
\frac{
(N-2)^{2}}{4}\left(1-\frac{1}{\tilde{S}^{\frac{N}{N-1}}}\right)
,
\frac{(N-2)^{2}}{4}
\right)$,
we know
\begin{equation}\label{34}
\begin{aligned}
S=\tilde{S}
\left(
1-\frac{4}{(N-2)^{2}}\zeta
\right)^{\frac{N-1}{N}}<1.
\end{aligned}
\end{equation}
Applying
$\zeta\in
\left(
\frac{
(N-2)^{2}}{4}\left(1-\frac{1}{(k+1)^{\frac{N-2}{N-1}}}\right)
,
\frac{(N-2)^{2}}{4}
\right)$,
we have
\begin{equation*}
\begin{aligned}
k+1
<
\left(
1-\frac{4}{(N-2)^{2}}\zeta
\right)^{-\frac{N-1}{N-2}}
=
\left(
\frac{\tilde{S}}{
S}
\right)^{\frac{N}{N-2}},
\end{aligned}
\end{equation*}
which gives that
\begin{equation}\label{35}
\begin{aligned}
\frac{\tilde{S} ^{\frac{N}{N+2-\alpha_{1}}}}{N}
(k+1)
^{-\frac{N-2}{N+2-\alpha_{1}}}
>
\frac{S ^{\frac{N}{N+2-\alpha_{1}}}}{N}.
\end{aligned}
\end{equation}
Inserting
(\ref{35})
into
(\ref{33}),
we know
\begin{equation}\label{36}
\begin{aligned}
\frac{\tilde{S} ^{\frac{N}{N+2-\alpha_{1}}}}{N}
(k+1)
^{-\frac{N-2}{N+2-\alpha_{1}}}
>
\frac{S ^{\frac{N}{N+2-\alpha_{1}}}}{N}
>
\frac{S ^{\frac{N}{2}}}{N}.
\end{aligned}
\end{equation}
According to
(\ref{36}),
Lemma \ref{lemma8}
and
Lemma \ref{lemma11},
we obtain
\begin{equation*}
\begin{aligned}
\tilde{c}_{0}
>
\frac{S ^{\frac{N}{2}}}{N}
>
c=c_{0}>0.
\end{aligned}
\end{equation*}

\noindent
{\bf Case
(3).}
Since
$\frac{2^{*}}{2}<2^{*}_{\alpha_{k}}<\cdots<2^{*}_{\alpha_{1}}$
and
$1
\leqslant
\tilde{S}_{\alpha_{j}}
=
\min
\left\{
\tilde{S}
,
\tilde{S}_{\alpha_{1}}
,\ldots,\tilde{S}_{\alpha_{k}}
\right\}$,
we have
\begin{equation}\label{37}
\begin{aligned}
\tilde{\kappa}
=
\frac{1}{\tilde{S} ^{\frac{2^{*} }{2}}}
+
\sum_{i=1}^{k}
\frac{1}
{\tilde{S}_{\alpha_{i}} ^{2^{*}_{\alpha_{i}}}}
\leqslant
\frac{k+1}{\tilde{S}_{\alpha_{j}} ^{\frac{2^{*} }{2}}}
<
k+1.
\end{aligned}
\end{equation}
Combining
(\ref{28}),
(\ref{37}),
$\tilde{\kappa}\geqslant1$
and
$-\frac{N-2}{N+2-\alpha_{1}}>-\frac{N-2}{2}$,
we get
\begin{equation}\label{38}
\begin{aligned}
\tilde{c}_{0}
\geqslant
\frac{1}{N}
\min
\left\{
\tilde{\kappa}
^{-\frac{N-2}{N+2-\alpha_{1}}}
,
\tilde{\kappa}
^{-\frac{N-2}{2}}
\right\}
=
\frac{1}{N}
\tilde{\kappa}
^{-\frac{N-2}{2}}
>
\frac{1}{N}
(k+1)^{-\frac{N-2}{2}}.
\end{aligned}
\end{equation}
Since
$\zeta\in
\left(
\frac{
(N-2)^{2}}{4}
\left(
1-
\frac{1}
{
(k+1)^{\frac{N-2}{N-1}}
\tilde{S}^{\frac{N}{N-1}}
}
\right)
,
\frac{(N-2)^{2}}{4}
\right)$,
we know
\begin{equation*}
\begin{aligned}
k+1
<&
\left(
\frac{1}{
\tilde{S}
\left(
1-\frac{4}{(N-2)^{2}}\zeta
\right)^{\frac{N-1}{N}}}
\right)^{\frac{N}{N-2}}
=
\left(
\frac{1}{
S}
\right)^{\frac{N}{N-2}},
\end{aligned}
\end{equation*}
which implies that
\begin{equation}\label{39}
\begin{aligned}
(k+1)^{-\frac{N-2}{2}}
>
S^{\frac{N}{2}}.
\end{aligned}
\end{equation}
According to
(\ref{38}),
(\ref{39}),
Lemma \ref{lemma8}
and
Lemma \ref{lemma11},
we obtain
\begin{equation*}
\begin{aligned}
\tilde{c}_{0}
>
\frac{S ^{\frac{N}{2}}}{N}
>
c=c_{0}>0.
\end{aligned}
\end{equation*}

\noindent
{\bf Case
(4).}
Similar to Case (3),
we have
$$\tilde{\kappa}< k+1,$$
which gives that
$$\{\tilde{\kappa}\leqslant 1\}\cap\{\tilde{\kappa}< k+1\}=\{\tilde{\kappa}\leqslant 1\}.$$
Combining
(\ref{28}),
$\tilde{\kappa}\leqslant 1$
and
$-\frac{N-2}{N+2-\alpha_{1}}>-\frac{N-2}{2}$,
we get
\begin{equation}\label{40}
\begin{aligned}
\tilde{c}_{0}
\geqslant
\frac{1}{N}
\min
\left\{
\tilde{\kappa}
^{-\frac{N-2}{N+2-\alpha_{1}}}
,
\tilde{\kappa}
^{-\frac{N-2}{2}}
\right\}
=
\frac{1}{N}
\tilde{\kappa}
^{-\frac{N-2}{N+2-\alpha_{1}}}
\geqslant
\frac{1}{N}.
\end{aligned}
\end{equation}
Since
$\zeta\in
\left(
\frac{
(N-2)^{2}}{4}\left(1-\frac{1}{\tilde{S}^{\frac{N}{N-1}}}\right)
,
\frac{(N-2)^{2}}{4}
\right)$,
we get
\begin{equation*}
\begin{aligned}
S=\tilde{S}
\left(
1-\frac{4}{(N-2)^{2}}\zeta
\right)^{\frac{N-1}{N}}<1,
\end{aligned}
\end{equation*}
which implies that
\begin{equation}\label{41}
\begin{aligned}
1
>
S
\Rightarrow
1
>
S^{\frac{N}{2}}.
\end{aligned}
\end{equation}
According to
(\ref{40}),
(\ref{41}),
Lemma \ref{lemma8}
and
Lemma \ref{lemma11},
we obtain
\begin{equation*}
\begin{aligned}
\tilde{c}_{0}
\geqslant
\frac{1}{N}
>
\frac{S ^{\frac{N}{2}}}{N}
>
c=c_{0}>0.
\end{aligned}
\end{equation*}
\qed
\section{The proof of Theorem \ref{theorem1}}
In this section,
we show the existence of nontrivial solution of problem
$(\mathcal{P})$.

\noindent
{\bf Proof of Theorem \ref{theorem1}:}
We divide our proof into five steps.

\noindent
{\bf Step 1.}
Since
$\{u_{n}\}$
is a bounded sequence in
$D^{1,2}(\mathbb{R}^{N})$,
up to a subsequence,
we assume that
\begin{align*}
&u_{n}\rightharpoonup u
~
\mathrm{in}
~
D^{1,2}(\mathbb{R}^{N}),~
u_{n}\rightarrow u
~
\mathrm{a.e. ~in}
~
\mathbb{R}^{N},\\
&u_{n}\rightarrow u
~
\mathrm{in}
~
L^{r}_{loc}(\mathbb{R}^{N})
~
\mathrm{for~all}
~
r\in[1,2^{*}).
\end{align*}
According to
Lemma
\ref{lemma4},
Lemma
\ref{lemma6}
and
Lemma
\ref{lemma12},
there exists
$C>0$
such that for any $n$ we get
$$
\|u_{n}\|_{\mathcal{L}^{2,N-2}(\mathbb{R}^{N})}\geqslant C>0.
$$
On the other hand,
since the sequence is bounded in
$D^{1,2}(\mathbb{R}^{N})$
and
$D^{1,2}(\mathbb{R}^{N})\hookrightarrow L^{2^{*}}(\mathbb{R}^{N})\hookrightarrow \mathcal{L}^{2,N-2}(\mathbb{R}^{N})$,
we have
$$
\|u_{n}\|_{\mathcal{L}^{2,N-2}(\mathbb{R}^{N})}\leqslant C,
$$
for some
$C>0$
independent of $n$.
Hence, there exists a positive constant which we denote again by $C$ such
that for any $n$ we obtain
$$
C
\leqslant
\|u_{n}\|_{\mathcal{L}^{2,N-2}(\mathbb{R}^{N})}\leqslant C^{-1}.
$$
Combining
the definition of Morrey space and above inequalities,
we deduce that for any
$n\in \mathbb{N}$
there exist
$\sigma_{n} > 0$
and
$x_{n}\in \mathbb{R}^{N}$
such that
$$
\frac{1}{\sigma_{n}^{2}}
\int_{B(x_{n},\sigma_{n})}
|u_{n}(y)|^{2}
\mathrm{d}y
\geqslant
\|u_{n}\|_{\mathcal{L}^{2,N-2}(\mathbb{R}^{N})}^{2}
-
\frac{C}{2n}
\geqslant
C_{4}>0.
$$
Let
$v_{n}(x)=\sigma_{n}^{\frac{N-2}{2}}u_{n}(x_{n}+\sigma_{n}x)$.
We may readily verify that
$$\tilde{I}_{\zeta}(v_{n})=I(u_{n})\rightarrow c,~
 \tilde{I}^{'}_{\zeta}(v_{n})\rightarrow0
~\mathrm{as}~n\rightarrow\infty,$$
where
\begin{equation*}
\begin{aligned}
\tilde{I}_{\zeta}(v_{n})
=&
\frac{1}{2}
\|v_{n}\|_{D}^{2}
-
\frac{\zeta}{2}
\int_{\mathbb{R}^{N}}
\frac{|v_{n}|^{2}}
{|x+\frac{x_{n}}{\sigma_{n}}|^{2}}
\mathrm{d}x\\
&-
\sum_{i=1}^{k}
\frac{1}{2\cdot2^{*}_{\alpha_{i}}}
\int_{\mathbb{R}^{N}}
\int_{\mathbb{R}^{N}}
\frac{|v_{n}(x)|^{2^{*}_{\alpha_{i}}}|v_{n}(y)|^{2^{*}_{\alpha_{i}}}}{|x-y|^{\alpha_{i}}}
\mathrm{d}x
\mathrm{d}y
-
\frac{1}{2^{*}}
\int_{\mathbb{R}^{N}}
|v_{n}|^{2^{*}}
\mathrm{d}x.
\end{aligned}
\end{equation*}
Now,
for all
$\varphi\in D^{1,2}(\mathbb{R}^{N})$,
we obtain
\begin{equation*}
\begin{aligned}
|\langle \tilde{I}^{'}_{\zeta}(v_{n}),\varphi\rangle|
=&
|\langle I^{'}(u_{n}),\bar{\varphi}\rangle|\\
\leqslant&
\|I^{'}(u_{n})\|_{D^{-1}}
\|\bar{\varphi}\|_{D}\\
=&
o(1)
\|\bar{\varphi}\|_{D},
\end{aligned}
\end{equation*}
where
$\bar{\varphi}=\sigma_{n}^{-\frac{N-2}{2}}\varphi(\frac{x-x_{n}}{\sigma_{n}})$.
Since
$\|\bar{\varphi}\|_{D}=\|\varphi\|_{D}$,
we get
$$\tilde{I}^{'}_{\zeta}(v_{n})\rightarrow0~
\mathrm{as}~
n\rightarrow\infty.$$
Thus there exists
$v$
such that
\begin{align*}
&v_{n}\rightharpoonup v
~
\mathrm{in}
~
D^{1,2}(\mathbb{R}^{N}),~
v_{n}\rightarrow v
~
\mathrm{a.e. ~in}
~
\mathbb{R}^{N},\\
&v_{n}\rightarrow v
~
\mathrm{in}
~
L^{r}_{loc}(\mathbb{R}^{N})~~\mathrm{for~all~}r\in[1,2^{*} ),
\end{align*}
and
$$
\int_{B(0,1)}
|v_{n}(y)|^{2}
\mathrm{d}y
=
\frac{1}{\sigma_{n}^{2}}
\int_{B(x_{n},\sigma_{n})}
|u_{n}(y)|^{2}
\mathrm{d}y
\geqslant
C_{4}>0.
$$
Hence,
$v\not\equiv0$.

\noindent
{\bf Step 2.}
Now,
we claim that
$\{\frac{x_{n}}{\sigma_{n}}\}$
is bounded.
If
$\frac{x_{n}}{\sigma_{n}}\rightarrow\infty$,
then for any
$\varphi\in D^{1,2}(\mathbb{R}^{N})$,
we get
\begin{equation}\label{42}
\begin{aligned}
\lim_{n\rightarrow\infty}
\int_{\mathbb{R}^{N}}
\frac{v_{n}\varphi}
{|x+\frac{x_{n}}{\sigma_{n}}|^{2}}
\mathrm{d}x
\rightarrow0.
\end{aligned}
\end{equation}
We  will show that
\begin{equation}\label{43}
\begin{aligned}
\langle \tilde{I}^{'}_{0}(v),\varphi\rangle=0.
\end{aligned}
\end{equation}
Since
$v_{n}\rightharpoonup v$
weakly in
$D^{1,2}(\mathbb{R}^{N})$,
we know
\begin{align}\label{44}
\lim_{n\rightarrow\infty}
\int_{\mathbb{R}^{N}}
\nabla v_{n}
\nabla \varphi
\mathrm{d}x
=\int_{\mathbb{R}^{N}}\nabla v
\nabla \varphi\mathrm{d}x.
\end{align}
By the Hardy--Littlewood--Sobolev inequality,
the Riesz potential defines a linear continuous map from
$L^{\frac{2N}{2N-\alpha_{i}}}(\mathbb{R}^{N})$
to
$L^{\frac{2N}{\alpha_{i}}}(\mathbb{R}^{N})$.
Since
$|v_{n}|^{2^{*}_{\alpha_{i}}}\rightharpoonup|v|^{2^{*}_{\alpha_{i}}}$
weakly in
$L^{\frac{2^{*}}{2^{*}_{\alpha_{i}}}}(\mathbb{R}^{N})$,
we have
\begin{equation}\label{45}
\int_{\mathbb{R}^{N}}
\frac{|v_{n}(y)|^{2^{*}_{\alpha_{i}}}}{|x-y|^{\alpha_{i}}}
\mathrm{d}y
\rightharpoonup
\int_{\mathbb{R}^{N}}
\frac{|v(y)|^{2^{*}_{\alpha_{i}}}}{|x-y|^{\alpha_{i}}}
\mathrm{d}y
~\mathrm{weakly~in}~L^{\frac{2N}{\alpha_{i}}}(\mathbb{R}^{N}).
\end{equation}
Now,
we show that
$|v_{n}|^{2^{*}_{\alpha_{i}}-2}
v_{n}
\varphi
\rightarrow
|v|^{2^{*}_{\alpha_{i}}-2}
v
\varphi$
in
$L^{\frac{2N}{2N-\alpha_{i}}}(\mathbb{R}^{N})$.
For any $\varepsilon>0$,
there exists $R>0$ large enough such that
\begin{equation}\label{46}
\begin{aligned}
&
\lim_{n\rightarrow\infty}
\int_{|x| >R}
\left|
|v_{n}|^{2^{*}_{\alpha_{i}}-2}
v_{n}
\varphi
\right|^{\frac{2N}{2N-\alpha_{i}}}
-
\left|
|v|^{2^{*}_{\alpha_{i}}-2}
v
\varphi
\right|^{\frac{2N}{2N-\alpha_{i}}}
\mathrm{d}x\\
\leqslant&
\lim_{n\rightarrow\infty}
\int_{|x| >R}
\left|
v_{n}
\right|^{(2^{*}_{\alpha_{i}}-1)\cdot\frac{2^{*}}{2^{*}_{\alpha_{i}}}}
\left|
\varphi
\right|^{\frac{2^{*}}{2^{*}_{\alpha_{i}}}}
\mathrm{d}x
+
\int_{|x| >R}
\left|
v
\right|^{(2^{*}_{\alpha_{i}}-1)\cdot\frac{2^{*}}{2^{*}_{\alpha_{i}}}}
\left|
\varphi
\right|^{\frac{2^{*}}{2^{*}_{\alpha_{i}}}}
\mathrm{d}x\\
\leqslant&
\lim_{n\rightarrow\infty}
\left(
\int_{|x| >R}
\left|
v_{n}
\right|^{2^{*}}
\mathrm{d}x
\right)^{1-\frac{1}{2^{*}_{\alpha_{i}}}}
\left(
\int_{|x| >R}
\left|
\varphi
\right|^{2^{*}}
\mathrm{d}x
\right)^{\frac{1}{2^{*}_{\alpha_{i}}}}\\
&+
\left(
\int_{|x| >R}
\left|
v
\right|^{2^{*}}
\mathrm{d}x
\right)^{1-\frac{1}{2^{*}_{\alpha_{i}}}}
\left(
\int_{|x| >R}
\left|
\varphi
\right|^{2^{*}}
\mathrm{d}x
\right)^{\frac{1}{2^{*}_{\alpha_{i}}}}\\
<&
\frac{\varepsilon}{2}.
\end{aligned}
\end{equation}
On the other hand,
by the boundedness of
$\{v_{n}\}$,
one has
\begin{equation*}
\begin{aligned}
\left(
\int_{|x|\leqslant R}
\left|
v_{n}
\right|^{2^{*}}
\mathrm{d}x
\right)^{1-\frac{1}{2^{*}_{\alpha_{i}}}}
\leqslant&
M.
\end{aligned}
\end{equation*}
where $M>0$ is a constant.
Let
$\Omega=\{x\in\mathbb{R}^{N}||x|\leqslant R\}$.
For any $\tilde{\varepsilon}>0$
there exists
$\delta>0$,
when
$E\subset\Omega$
with
$|E|<\delta$.
We obtain
\begin{align*}
\int_{E}
\left|
|v_{n}|^{2^{*}_{\alpha_{i}}-2}
v_{n}
\varphi
\right|^{\frac{2N}{2N-\alpha_{i}}}
\mathrm{d}x
=&
\int_{E}
\left|
v_{n}
\right|^{(2^{*}_{\alpha_{i}}-1)\cdot\frac{2^{*}}{2^{*}_{\alpha_{i}}}}
\left|
\varphi
\right|^{\frac{2^{*}}{2^{*}_{\alpha_{i}}}}
\mathrm{d}x\\
\leqslant&
\left(
\int_{E}
\left|
v_{n}
\right|^{2^{*}}
\mathrm{d}x
\right)^{1-\frac{1}{2^{*}_{\alpha_{i}}}}
\left(
\int_{E}
\left|
\varphi
\right|^{2^{*}}
\mathrm{d}x
\right)^{\frac{1}{2^{*}_{\alpha_{i}}}}\\
<&M\tilde{\varepsilon}.
\end{align*}
where the last inequality is from the absolutely continuity of
$\int_{E}
\left|
\varphi
\right|^{2^{*}}
\mathrm{d}x$.
Moreover,
$|v_{n}|^{2^{*}_{\alpha_{i}}-2}
v_{n}
\varphi
\rightarrow
|v|^{2^{*}_{\alpha_{i}}-2}
v
\varphi$
a.e. in $\mathbb{R}^{N}$ as $n\rightarrow\infty$.
Thus,
by the Vitali convergence Theorem,
we get
\begin{equation}\label{47}
\begin{aligned}
\lim_{n\rightarrow\infty}
\int_{|x|\leqslant R}
\left|
|v_{n}|^{2^{*}_{\alpha_{i}}-2}
v_{n}
\varphi
\right|^{\frac{2N}{2N-\alpha_{i}}}
\mathrm{d}x
=
\int_{|x|\leqslant R}
\left|
|v|^{2^{*}_{\alpha_{i}}-2}
v
\varphi
\right|^{\frac{2N}{2N-\alpha_{i}}}
\mathrm{d}x.
\end{aligned}
\end{equation}
It follows from
\eqref{46}
and
\eqref{47}
that
\begin{equation*}
\begin{aligned}
&
\lim_{n\rightarrow\infty}
\left|
\int_{\mathbb{R}^{N}}
\left|
|v_{n}|^{2^{*}_{\alpha_{i}}-2}
v_{n}
\varphi
\right|^{\frac{2N}{2N-\alpha_{i}}}
-
\left|
|v|^{2^{*}_{\alpha_{i}}-2}
v
\varphi
\right|^{\frac{2N}{2N-\alpha_{i}}}
\mathrm{d}x
\right|\\
\leqslant&
\lim_{n\rightarrow\infty}
\left|
\int_{|x|\leqslant R}
\left|
|v_{n}|^{2^{*}_{\alpha_{i}}-2}
v_{n}
\varphi
\right|^{\frac{2N}{2N-\alpha_{i}}}
-
\left|
|v|^{2^{*}_{\alpha_{i}}-2}
v
\varphi
\right|^{\frac{2N}{2N-\alpha_{i}}}
\mathrm{d}x
\right|\\
&+
\lim_{n\rightarrow\infty}
\left|
\int_{|x|> R}
\left|
|v_{n}|^{2^{*}_{\alpha_{i}}-2}
v_{n}
\varphi
\right|^{\frac{2N}{2N-\alpha_{i}}}
-
\left|
|v|^{2^{*}_{\alpha_{i}}-2}
v
\varphi
\right|^{\frac{2N}{2N-\alpha_{i}}}
\mathrm{d}x
\right|\\
<&
\varepsilon.
\end{aligned}
\end{equation*}
This implies that
\begin{equation}\label{48}
\begin{aligned}
\lim_{n\rightarrow\infty}
\int_{\mathbb{R}^{N}}
\left|
|v_{n}|^{2^{*}_{\alpha_{i}}-2}
v_{n}
\varphi
\right|^{\frac{2N}{2N-\alpha_{i}}}
\mathrm{d}x
=
\int_{\mathbb{R}^{N}}
\left|
|v|^{2^{*}_{\alpha_{i}}-2}
v
\varphi
\right|^{\frac{2N}{2N-\alpha_{i}}}
\mathrm{d}x
\end{aligned}
\end{equation}
Combining
\eqref{45}
and
\eqref{48},
we have
\begin{equation}\label{49}
\begin{aligned}
&\lim_{n\rightarrow\infty}
\int_{\mathbb{R}^{N}}
\int_{\mathbb{R}^{N}}
\frac{
|v_{n}(y)|^{2^{*}_{\alpha_{i}}}
|v_{n}(x)|^{2^{*}_{\alpha_{i}}-2}
v_{n}(x)
\varphi(x)
}{|x-y|^{\alpha_{i}}}
\mathrm{d}y
\mathrm{d}x\\
=&
\int_{\mathbb{R}^{N}}
\int_{\mathbb{R}^{N}}
\frac{
|v(y)|^{2^{*}_{\alpha_{i}}}
|v(x)|^{2^{*}_{\alpha_{i}}-2}
v(x)
\varphi(x)
}{|x-y|^{\alpha_{i}}}
\mathrm{d}y
\mathrm{d}x.
\end{aligned}
\end{equation}
Similarly,
we get
\begin{equation}\label{49'}
\begin{aligned}
\lim_{n\rightarrow\infty}
\int_{\mathbb{R}^{N}}
|v_{n}|^{2^{*}-2}
v_{n}
\varphi
\mathrm{d}x
=
\int_{\mathbb{R}^{N}}
|v|^{2^{*}-2}
v
\varphi
\mathrm{d}x.
\end{aligned}
\end{equation}
Applying
$\lim\limits_{n\rightarrow\infty}\langle \tilde{I}_{\zeta}^{'}(v_{n}),\varphi\rangle\rightarrow0$,
\eqref{42},
\eqref{44},
\eqref{49}
and
\eqref{49'}
we know
\begin{equation}\label{50}
\begin{aligned}
\langle \tilde{I}^{'}_{0}(v),\varphi\rangle=0.
\end{aligned}
\end{equation}
Moreover,
according to \eqref{50}
and
$v\not\equiv0$,
we get that
$$v\in \widetilde{\mathcal{N}}.$$
By Br\'{e}zis--Lieb lemma \cite[Lemma 2.2]{Gao2016},
we have
\begin{equation*}
\begin{aligned}
&\int_{\mathbb{R}^{N}}
\int_{\mathbb{R}^{N}}
\frac{|v_{n}(x)|^{2^{*}_{\alpha}}|v_{n}(y)|^{2^{*}_{\alpha}}}{|x-y|^{\alpha}}
\mathrm{d}x
\mathrm{d}y
-
\int_{\mathbb{R}^{N}}
\int_{\mathbb{R}^{N}}
\frac{|v_{n}(x)-v(x)|^{2^{*}_{\alpha}}|v_{n}(y)-v(y)|^{2^{*}_{\alpha}}}{|x-y|^{\alpha}}
\mathrm{d}x
\mathrm{d}y\\
=&
\int_{\mathbb{R}^{N}}
\int_{\mathbb{R}^{N}}
\frac{|v(x)|^{2^{*}_{\alpha}}|v(y)|^{2^{*}_{\alpha}}}{|x-y|^{\alpha}}
\mathrm{d}x
\mathrm{d}y
+
o(1),
\end{aligned}
\end{equation*}
which implies that
\begin{equation}\label{51}
\begin{aligned}
\int_{\mathbb{R}^{N}}
\int_{\mathbb{R}^{N}}
\frac{|v_{n}(x)|^{2^{*}_{\alpha}}|v_{n}(y)|^{2^{*}_{\alpha}}}{|x-y|^{\alpha}}
\mathrm{d}x
\mathrm{d}y
\geqslant
\int_{\mathbb{R}^{N}}
\int_{\mathbb{R}^{N}}
\frac{|v(x)|^{2^{*}_{\alpha}}|v(y)|^{2^{*}_{\alpha}}}{|x-y|^{\alpha}}
\mathrm{d}x
\mathrm{d}y
+
o(1).
\end{aligned}
\end{equation}
Similarly,
we get
\begin{equation}\label{52}
\begin{aligned}
\int_{\mathbb{R}^{N}}
|v_{n}|^{2^{*}}
\mathrm{d}x
\geqslant
\int_{\mathbb{R}^{N}}
|v|^{2^{*}}
\mathrm{d}x
+
o(1).
\end{aligned}
\end{equation}
Set
\begin{align*}
K(u)
=
\sum_{i=1}^{k}
\left(
\frac{1}{2}
-
\frac{1}{2\cdot2^{*}_{\alpha_{i}}}
\right)
\int_{\mathbb{R}^{N}}
\int_{\mathbb{R}^{N}}
\frac{|u(x)|^{2^{*}_{\alpha_{i}}}|u(y)|^{2^{*}_{\alpha_{i}}}}{|x-y|^{\alpha_{i}}}
\mathrm{d}x
\mathrm{d}y
+
\left(
\frac{1}{2}
-
\frac{1}{2^{*}}
\right)
\int_{\mathbb{R}^{N}}
|u|^{2^{*}}
\mathrm{d}x.
\end{align*}
Applying
Lemma \ref{lemma13},
Lemma \ref{lemma8},
(\ref{51}),
(\ref{52}),
$v\in \widetilde{\mathcal{N}}$
and
Lemma \ref{lemma9},
we obtain
\begin{align*}
\tilde{c}_{0}
>
c_{0}
=c
=&
I(v_{n})
-
\frac{1}{2}
\langle I^{'}(v_{n}),v_{n}\rangle\\
=&
\lim_{n\rightarrow\infty}
K(v_{n})
+o(1)\\
\geqslant&
K(v)
+o(1)\\
=&
\tilde{I}_{0}(v)
-
\frac{1}{2}
\langle \tilde{I}_{0}^{'}(v),v\rangle
=
 \tilde{I}_{0}(v)
\geqslant
\tilde{c}_{0},
\end{align*}
which yields a contradiction.
Hence, $\{\frac{x_{n}}{\sigma_{n}}\}$ is bounded.

\noindent
{\bf Step 3.}
In this step,
we study another $(PS)_{c}$ sequence of $I$.
Let
$\tilde{v}_{n}(x)=\sigma_{n}^{\frac{N-2}{2}}u_{n}(\sigma_{n}x)$.
Then we can
verify that
$$I(\tilde{v}_{n})=I(u_{n})\rightarrow c,~
I^{'}(\tilde{v}_{n})\rightarrow 0
~\mathrm{as}~n\rightarrow\infty.$$
Arguing as before, we have
\begin{align*}
&\tilde{v}_{n}\rightharpoonup \tilde{v}
~
\mathrm{in}
~
D^{1,2}(\mathbb{R}^{N}),~
\tilde{v}_{n}\rightarrow \tilde{v}
~
\mathrm{a.e. ~in}
~
\mathbb{R}^{N},\\
&\tilde{v}_{n}\rightarrow \tilde{v}
~
\mathrm{in}
~
L^{r}_{loc}(\mathbb{R}^{N})~~\mathrm{for~all~}r\in[1,2^{*}).
\end{align*}
Since $\{\frac{x_{n}}{\sigma_{n}}\}$ is bounded,
there exists
$\tilde{R}>0$ such that
$$
\int_{B(0,\tilde{R})}
|\tilde{v}_{n}(y)|^{2}
\mathrm{d}y
>
\int_{B(\frac{x_{n}}{\sigma_{n}},1)}
|\tilde{v}_{n}(y)|^{2}
\mathrm{d}y
=
\frac{1}{\sigma_{n}^{2}}
\int_{B(x_{n},\sigma_{n})}
|u_{n}(y)|^{2}
\mathrm{d}y
\geqslant
C_{4}>0.$$
As a result,
$\tilde{v}\not\equiv0$.

\noindent
{\bf Step 4.}
In this step,
we show
$\tilde{v}_{n}\rightarrow \tilde{v}$
strongly in
$D^{1,2}(\mathbb{R}^{N})$.
Similar to Step 2,
we know that
\begin{equation}\label{53}
\begin{aligned}
\langle I^{'}(\tilde{v}),\varphi\rangle=0.
\end{aligned}
\end{equation}
Applying
\eqref{51}
--
\eqref{53}
and
Lemma \ref{lemma8},
we obtain
\begin{equation}\label{54}
\begin{aligned}
c
=&
I(\tilde{v}_{n})
-
\frac{1}{2}
\langle I^{'}(\tilde{v}_{n}),\tilde{v}_{n}\rangle
+
o(1)\\
=&
\lim_{n\rightarrow\infty}
K(v_{n})
+o(1)\\
\geqslant&
\lim_{n\rightarrow\infty}
K(v)
+o(1)\\
=&
I(\tilde{v})
-
\frac{1}{2}
\langle I^{'}(\tilde{v}),\tilde{v}\rangle
=
I(\tilde{v})
\geqslant
c.
\end{aligned}
\end{equation}
\textcolor{red}{Therefore, the inequalities above have to be equalities}.
We know
\begin{align*}
\lim\limits_{n\rightarrow\infty}
K(\tilde{v}_{n})
=
K(\tilde{v}).
\end{align*}
By using Br\'{e}zis--Lieb lemma again,
we have
\begin{align*}
\lim\limits_{n\rightarrow\infty}
K(\tilde{v}_{n})
-
\lim\limits_{n\rightarrow\infty}
K(\tilde{v}_{n}-\tilde{v})
=
K(\tilde{v})+o(1).
\end{align*}
Hence,
we conclude that
\begin{align*}
\lim\limits_{n\rightarrow\infty}
K(\tilde{v}_{n}-\tilde{v})
=
0,
\end{align*}
which implies that
\begin{equation}\label{55}
\begin{aligned}
&
\lim\limits_{n\rightarrow\infty}
\int_{\mathbb{R}^{N}}
\int_{\mathbb{R}^{N}}
\frac{|\tilde{v}_{n}(x)-\tilde{v}(x)|^{2^{*}_{\alpha_{i}}}|\tilde{v}_{n}(y)-\tilde{v}(y)|^{2^{*}_{\alpha_{i}}}}{|x-y|^{\alpha_{i}}}
\mathrm{d}x
\mathrm{d}y
=0
~\mathrm{for~all}~
i=1,\ldots,k,\\
&
\lim\limits_{n\rightarrow\infty}
\int_{\mathbb{R}^{N}}
|\tilde{v}_{n}(x)-\tilde{v}(x)|^{2^{*}}
\mathrm{d}x
=0.
\end{aligned}
\end{equation}
According to
$\langle I^{'}(\tilde{v}_{n}),\tilde{v}_{n}\rangle=o(1)$,
$\langle I^{'}(\tilde{v}),\tilde{v}\rangle=0$
and Br\'{e}zis--Lieb lemma,
we obtain
\begin{equation*}
\begin{aligned}
o(1)=&
\langle I^{'}(\tilde{v}_{n}),\tilde{v}_{n}\rangle
-
\langle I^{'}(\tilde{v}),\tilde{v}\rangle\\
=&
\|\tilde{v}_{n}-\tilde{v}\|_{\zeta}^{2}
-
\sum_{i=1}^{k}
\int_{\mathbb{R}^{N}}
\int_{\mathbb{R}^{N}}
\frac{|\tilde{v}_{n}(x)-\tilde{v}(x)|^{2^{*}_{\alpha_{i}}}|\tilde{v}_{n}(y)-\tilde{v}(y)|^{2^{*}_{\alpha_{i}}}}{|x-y|^{\alpha_{i}}}
\mathrm{d}x
\mathrm{d}y\\
&-
\int_{\mathbb{R}^{N}}
|\tilde{v}_{n}-\tilde{v}|^{2^{*}}
\mathrm{d}x+o(1),
\end{aligned}
\end{equation*}
which implies that
\begin{equation}\label{56}
\begin{aligned}
&
\lim\limits_{n\rightarrow\infty}
\|\tilde{v}_{n}-\tilde{v}\|_{\zeta}^{2}\\
=&
\lim\limits_{n\rightarrow\infty}
\sum_{i=1}^{k}
\int_{\mathbb{R}^{N}}
\int_{\mathbb{R}^{N}}
\frac{|\tilde{v}_{n}(x)-\tilde{v}(x)|^{2^{*}_{\alpha_{i}}}|\tilde{v}_{n}(y)-\tilde{v}(y)|^{2^{*}_{\alpha_{i}}}}{|x-y|^{\alpha_{i}}}
\mathrm{d}x
\mathrm{d}y\\
&+
\lim\limits_{n\rightarrow\infty}
\int_{\mathbb{R}^{N}}
|\tilde{v}_{n}-\tilde{v}|^{2^{*}}
\mathrm{d}x+o(1).
\end{aligned}
\end{equation}
Combining
(\ref{55})
and
(\ref{56}),
we get
\begin{equation*}
\begin{aligned}
\lim\limits_{n\rightarrow\infty}
\|\tilde{v}_{n}-\tilde{v}\|_{\zeta}^{2}
\rightarrow
o(1).
\end{aligned}
\end{equation*}
Since
$\tilde{v}\not\equiv0$,
we know that
$\tilde{v}_{n}\rightarrow \tilde{v}$
strongly in
$D^{1,2}(\mathbb{R}^{N})$.

\noindent
{\bf Step 5.}
Here we show the properties of the solution.
By using \eqref{54} again,
we know
$I(\tilde{v})=c$,
which means that $\tilde{v}$ is a nontrivial solution of problem $(\mathcal{P})$ at the
energy level $c$.
Since
$I$ is even,
we have
$$c=I(\tilde{v})=I(|\tilde{v}|)~\mathrm{and}~\langle I^{'}(\tilde{v}),\tilde{v}\rangle=\langle I^{'}(|\tilde{v}|),|\tilde{v}|\rangle=0.$$
Then
$|\tilde{v}|$
is also a
critical point of
$I$.
Hence,
we can choose
$\tilde{v}\geqslant0$.
By the Kelvin transformation,
we have
\begin{equation}\label{57}
\begin{aligned}
\tilde{\tilde{v}}(x)
=
\frac{1}{|x|^{N-2}}
\tilde{v}
\left(
\frac{x}{|x|^{2}}
\right).
\end{aligned}
\end{equation}
It is well known that
\begin{equation}\label{58}
\begin{aligned}
-
\Delta
\tilde{\tilde{v}}(x)
=
\frac{1}{|x|^{N+2}}
(-
\Delta
\tilde{v})
\left(
\frac{x}{|x|^{2}}
\right),
\end{aligned}
\end{equation}
and
\begin{equation}\label{59}
\begin{aligned}
\frac{\tilde{\tilde{v}}(x)}{|x|^{2}}
=
\frac{1}{|x|^{N+2}}
\frac{
\tilde{v}
\left(
\frac{x}{|x|^{2}}
\right)
}
{\left|\frac{x}{|x|^{2}}\right|}.
\end{aligned}
\end{equation}
The following identity is very useful.
For
$\forall x,y\in \mathbb{R}^{N}\backslash\{0\}$,
we get
\begin{equation}\label{60}
\begin{aligned}
\frac{1}
{
\left|
\frac{x}{|x|^{2}}-\frac{y}{|y|^{2}}
\right|^{\alpha_{i}}}
\cdot
\frac{1}{|xy|^{\alpha_{i}}}
=&
\frac{1}
{
\left|
\frac{x\cdot y^{2}-y\cdot x^{2}}{(xy)^{2}}
\right|^{\alpha_{i}}}
\cdot
\frac{1}{|xy|^{\alpha_{i}}}\\
=&
\frac{1}
{
\left|
\frac{x\cdot y^{2}-y\cdot x^{2}}{xy}
\right|^{\alpha_{i}}}\\
=&
\frac{1}
{
\left|
y-x
\right|^{\alpha_{i}}}.
\end{aligned}
\end{equation}
Set
$z=\frac{y}{|y|^{2}}$.
Applying
(\ref{57})
and
(\ref{60}),
we have
\begin{equation}\label{61}
\begin{aligned}
\int_{\mathbb{R}^{N}}
\frac{|\tilde{\tilde{v}}(y)|^{\frac{2N-\alpha_{i}}{N-2}}}
{|x-y|^{\alpha_{i}}}
\mathrm{d}y
=&
\int_{\mathbb{R}^{N}}
\frac{|\tilde{v}
\left(
\frac{y}{|y|^{2}}
\right)|^{\frac{2N-\alpha_{i}}{N-2}}}
{|x-y|^{\alpha_{i}}}
\cdot
\frac{1}{|y|^{2N-\alpha_{i}}}
\mathrm{d}y
~(\mathrm{by}~(\ref{57}))\\
=&
\int_{\mathbb{R}^{N}}
\frac{|\tilde{v}
\left(
\frac{y}{|y|^{2}}
\right)|^{\frac{2N-\alpha_{i}}{N-2}}}
{|\frac{x}{|x|^{2}}-\frac{y}{|y|^{2}}|^{\alpha_{i}}}
\cdot
\frac{1}{|xy|^{\alpha_{i}}\cdot|y|^{2N-\alpha_{i}}}
\mathrm{d}y
(\mathrm{by}~(\ref{60}))\\
=&
\frac{1}{|x|^{\alpha_{i}}}
\int_{\mathbb{R}^{N}}
\frac{|\tilde{v}
\left(
\frac{y}{|y|^{2}}
\right)|^{\frac{2N-\alpha_{i}}{N-2}}}
{|\frac{x}{|x|^{2}}-\frac{y}{|y|^{2}}|^{\alpha_{i}}}
\cdot
\frac{1}{|y|^{2N}}
\mathrm{d}y\\
=&
\frac{1}{|x|^{\alpha_{i}}}
\int_{\mathbb{R}^{N}}
\frac{|\tilde{v}
\left(
z
\right)|^{\frac{2N-\alpha_{i}}{N-2}}}
{|\frac{x}{|x|^{2}}-z|^{\alpha_{i}}}
\mathrm{d}z
~(\mathrm{set}~z=\frac{y}{|y|^{2}}).
\end{aligned}
\end{equation}
Therefore,
by using
\eqref{58},
\eqref{59}
and
\eqref{61},
we get
\begin{equation*}
\begin{aligned}
&
\frac{1}{|x|^{N+2}}
(-
\Delta
\tilde{v})
\left(
\frac{x}{|x|^{2}}
\right)
-
\frac{\zeta}{|x|^{N+2}}
\cdot
\frac{
\tilde{v}
\left(
\frac{x}{|x|^{2}}
\right)
}
{\left|\frac{x}{|x|^{2}}\right|^{2}}\\
=&
\sum\limits^{k}_{i=1}
\frac{1}{|x|^{\alpha_{i}}}
\int_{\mathbb{R}^{N}}
\frac{|\tilde{v}
\left(
z
\right)|^{\frac{2N-\alpha_{i}}{N-2}}}
{|\frac{x}{|x|^{2}}-z|^{\alpha_{i}}}
\mathrm{d}z
\cdot
\left|
\frac{1}{|x|^{N-2}}
\tilde{v}
\left(
\frac{x}{|x|^{2}}
\right)
\right|^{\frac{4-\alpha_{i}}{N-2}}
\cdot
\frac{1}{|x|^{N-2}}
\tilde{v}
\left(
\frac{x}{|x|^{2}}
\right)\\
&+
\left|
\frac{1}{|x|^{N-2}}
\tilde{v}
\left(
\frac{x}{|x|^{2}}
\right)
\right|^{\frac{4}{N-2}}
\cdot
\frac{1}{|x|^{N-2}}
\tilde{v}
\left(
\frac{x}{|x|^{2}}
\right)\\
=&
\frac{1}{|x|^{N+2}}
\sum\limits^{k}_{i=1}
\int_{\mathbb{R}^{N}}
\frac{|\tilde{v}
\left(
z
\right)|^{\frac{2N-\alpha_{i}}{N-2}}}
{|\frac{x}{|x|^{2}}-z|^{\alpha_{i}}}
\mathrm{d}z
\cdot
\left|
\tilde{v}
\left(
\frac{x}{|x|^{2}}
\right)
\right|^{\frac{4-\alpha_{i}}{N-2}}
\cdot
\tilde{v}
\left(
\frac{x}{|x|^{2}}
\right)\\
&+
\frac{1}{|x|^{N+2}}
\left|
\tilde{v}
\left(
\frac{x}{|x|^{2}}
\right)
\right|^{\frac{4}{N-2}}
\cdot
\tilde{v}
\left(
\frac{x}{|x|^{2}}
\right),
\end{aligned}
\end{equation*}
which gives that
\begin{equation*}
\begin{aligned}
-
\Delta
\tilde{\tilde{v}}
-
\zeta
\frac{\tilde{\tilde{v}}}{|x|^{2}}
=&
\sum\limits^{k}_{i=1}
\left(
\int_{\mathbb{R}^{N}}
\frac{|\tilde{\tilde{v}}|^{2^{*}_{\alpha_{i}}}}
{|x-y|^{\alpha_{i}}}
\mathrm{d}y
\right)
\left|
\tilde{\tilde{v}}
\right|^{2^{*}_{\alpha_{i}}-2}
\tilde{\tilde{v}}
+
\left|
\tilde{\tilde{v}}
\right|^{2^{*}-2}
\tilde{\tilde{v}},
~\mathrm{in}~\mathbb{R}^{N}\backslash\{0\}.
\end{aligned}
\end{equation*}
\qed

\section*{Open Problem}
During the preparation of the manuscript we faced one problem which is worth to be tackled in forthcoming investigation.

We just study the case of $\tilde{S}_{\alpha}\geqslant1$ (see $(H_{2})$),
it is nature to ask the case of
$\tilde{S}_{\alpha}\in(0,1)$ (see Fig 1.).


\begin{thebibliography}{00}

\bibitem{O.Alves2017}
C. O. Alves, F. Gao, M. Squassina, M. Yang,
{\it Singularly Perturbed critical Choquard equations},
J. Differential Equations
{\bf 263} (2017) 3943--3988.

\bibitem{Bhakta2015}
M. Bhakta,
{\it A note on semilinear elliptic equation with biharmonic operator and multiple critical nonlinearities},
Adv. Nonlinear Stud. {\bf 15} (2015)  835--848.

\bibitem{Bellazzini2016}
J. Bellazzini, M. Ghimenti, C. Mercuri, V. Moroz, J. Van Schaftingen,
{\it Sharp Gagliardo--Nirenberg inequalities in fractional	Coulomb--Sobolev spaces},
T. Am. Math. Soc.
(2017);
Doi: https://doi.org/10.1090/tran/7426.

\bibitem{d'Avenia2015}
P. D'Avenia,
G. Siciliano,
M. Squassina,
{\it On fractional Choquard equations},
Math. Mod. and Meth. Appl. S.
{ \bf 25}  (2015) 1447--1476.

\bibitem{Pucci2009}
R. Filippucci, P. Pucci, F. Robert,
{\it On a p--Laplace equation with multiple critical nonlinearities},
J.  Math. Pures Appl.
{\bf 91} (2009) 156--177.

\bibitem{Felli2006}
V. Felli, S. Terracini,
{\it Elliptic Equations with Multi--Singular Inverse--Square Potentials and Critical Nonlinearity},
Commun. Part. Diff. Eq. {\bf 31} (2006) 469--495.

\bibitem{Gao2016}
F. Gao, M. Yang,
{\it On the Brezis--Nirenberg type critical problem for nonlinear Choquard equation},
Sci. China Math. (2016);	Doi: 10.1007/s11425--016--9067--5.

\bibitem{Gao2017JMAA}
F. Gao, M. Yang,
{\it On nonlocal Choquard equations with Hardy--Littlewood--Sobolev critical exponetns},
J. Math. Anal. Appl.
{\bf 448} (2017) 1006--1041.

\bibitem{Ghoussoub2016}
N. Ghoussoub, S. Shakerian,
{\it Borderline variational problems involving fractional Laplacians and critical singularities},
Adv. Nonlinear Stud.
{\bf 15} (2016) 527--555.

\bibitem{Lieb2001}
E. Lieb, M. Loss,
{\it Analysis, Gradute Studies in Mathematics},
AMS, Providence (2001).

\bibitem{Li2012}
Y. Li, C. S. Lin,
{\it A nonlinear elliptic pde with two sobolev--hardy critical exponents},
Arch. Ration. Mech. An. {\bf 203} (2012) 943--968.

\bibitem{Mukherjee2017Fractional}
T. Mukherjee, K. Sreenadh,
{\it Fractional Choquard equation with critical nonlinearities},
Nonlinear Differ. Equ. Appl.
{\bf 24} (2017) 63.

\bibitem{Mukherjee2017}
T. Mukherjee, K. Sreenadh,
{\it Positive solutions for nonlinear Choquard equation with singular nonlinearity},
Complex Var. Elliptic
{\bf 62} (2017) 1044--1071.

\bibitem{Mercuri2016}
C. Mercuri, V. Moroz, J. Van Schaftingen,
{\it Groundstates and radial solutions to nonlinear Schr\"{o}dinger--Poisson--Slater equations at the critical frequency},
Calc. Var.
{\bf 55} (2016) 146.

\bibitem{Moroz2016}
V. Moroz, J. Van Schaftingen,
{\it A guide to the Choquard equation},
J. Fix. Point Theory A.
(2016) 1--41.

\bibitem{Palatucci2014}
G. Palatucci, A. Pisante,
{\it Improved Sobolev embeddings, profile decomposition, and concentration--compactness for fractional Sobolev spaces},
Calc. Var.
{\bf 50} (2014) 799--829.

\bibitem{Pekar1954}
S. Pekar,
{\it Untersuchung\"{u}ber die Elektronentheorie der Kristalle},
Akademie Verlag, Berlin, (1954).

\bibitem{Penrose1996}
R. Penrose,
{\it On gravity's role in quantum sstate reduction},
Gen. Relativity Gravitation
{\bf 28} (1996) 581--600.

\bibitem{Su2018FractionalLaplacian}
Y. Su,
H. Chen,
{\it Finite many critical problems involving Fractional
Laplacians in $\mathbb{R}^{N}$},
subbmitted to Computers Mathematics with Applications.

\bibitem{Su2018pLaplacian}
Y. Su,
H. Chen,
{\it The minimizing problem involving p--Laplacian and Hardy--Littlewood--Sobolev upper critical exponent},
subbmitted to Electronic Journal of Qualitative Theory of Differential Equations.

\bibitem{Swanson1992}
C. A. Swanson,
{\it The best sobolev constant},
Appl. Anal. {\bf 47} (1992) 227--239.

\bibitem{Seok2018}
J. Seok,
{\it Nonlinear Choquard equations: Doubly critical cases},
Appl. Math. Lett.
{\bf 76} (2018) 148--156.

\bibitem{ZhangBL2017}
L. Wang,
B. Zhang,
H. Zhang,
{\it Fractional Laplacian system involving doubly critical nonlinearities in $\mathbb{R}^{N}$},
Electron. J. Qual. Theo.
{\bf 57} (2017) 1--17.

\bibitem{Willem1996}
M. Willem,
{\it Minimax theorems},
Birkh\"{a}user, Boston, (1996).

\bibitem{Yang2017}
J. Yang,
F. Wu,
{\it Doubly critical problems involving fractional laplacians in $\mathbb{R}^{N}$},
Adv. Nonlinear Stud.
(2017); Doi: https://doi.org/10.1515/ans--2016--6012.

\bibitem{Zhong2016}
X. Zhong, W. Zou,
{\it A perturbed nonlinear elliptic PDE with two Hardy--Sobolev critical exponents},
Commun. Contemp. Math. (2016) 1550061.
\end{thebibliography}
\end{document}